\documentclass[mathpazo]{cicp}
\usepackage{float}
\usepackage{graphicx,epstopdf}
\graphicspath{{./}{../}}
\usepackage{float}
\usepackage{algorithm}
\usepackage{multirow}
\usepackage{epstopdf}
\usepackage{amsmath}
\usepackage{booktabs}
\usepackage{makecell}
\usepackage{color}
\numberwithin{table}{section}
\numberwithin{figure}{section}

\begin{document}
\title{Kernel-learning parameter prediction and evaluation in algebraic multigrid method for several PDEs}



 \author[Zhang Z R et.~al.]{Junyue Luo\affil{}, Xiaoqiang Yue\affil{}, Fangfang Zhang\affil{} and Juan Zhang\affil{}\comma\corrauth}
 \address{\affilnum{}\ Hunan Key Laboratory for Computation and Simulation in Science and Engineering, Key  Laboratory of Intelligent Computing and Information Processing of Ministry of Education, School of Mathematics and Computational Science, Xiangtan University, Xiangtan, Hunan, China, 411105.}
 \emails{{\tt zhangjuan@xtu.edu.cn} (Juan Zhang)}
 \footnote{Funding. The National Key R\&D Program of China (2023YFB3001604).}

\begin{abstract}
This paper explores the application of kernel learning methods for parameter prediction and evaluation in the Algebraic Multigrid Method (AMG), focusing on several Partial Differential Equation (PDE) problems. AMG is an efficient iterative solver for large-scale sparse linear systems, particularly those derived from elliptic and parabolic PDE discretizations. However, its performance heavily relies on numerous parameters, which are often set empirically{ , and different parameter settings can significantly impact the effectiveness of AMG.}
 Traditional parameter optimization methods are either computationally expensive or lack theoretical support. To address this, we propose a Gaussian Process Regression (GPR)-based strategy to optimize AMG parameters and introduce evaluation metrics to assess their effectiveness. Trained on small-scale datasets, GPR predicts nearly optimal parameters, bypassing the time-consuming parameter sweeping process. We also use kernel learning techniques to build a kernel function library and determine the optimal kernel function through linear combination, enhancing prediction accuracy. In numerical experiments, we tested typical PDEs such as  Poisson equation, Parabolic equation, Diffusion equation, and Helmholtz equation. Results show that GPR-predicted parameters match grid search results in iteration steps while significantly reducing computational time. A comprehensive analysis using metrics like mean squared error, prediction interval coverage, and Bayesian information criterion confirms GPR's efficiency and reliability. These findings validate GPR's effectiveness in AMG parameter optimization and provide theoretical support for AMG's practical application.
\end{abstract}


\maketitle

   \section{Introduction}\label{sec1}

We consider the sparse linear system defined by
\begin{equation}\label{1}
	Au = b,
\end{equation}
where $A$ is a non-singular matrix and $b$ is a given vector. Solving this system is a cornerstone in scientific and engineering computing, particularly in the realm of elliptic partial differential equations (PDEs). These PDEs are notable for their ability to depict steady-state or equilibrium conditions. Notable instances include Laplace equation and Poisson equation, which are pivotal in modeling electric potential and gravitational fields, respectively. Numerically addressing such equations often leads to large-scale sparse linear systems, posing significant computational hurdles. Parabolic PDEs, which describe physical processes evolving over time, such as the heat conduction equation, also generate time-dependent linear systems. These systems necessitate the solution of sparse linear systems at each time step, similar to those arising in elliptic PDEs. The Helmholtz equation, a linear elliptic PDE with extensive applications in acoustics and electromagnetism, is particularly challenging due to its oscillatory solutions and the necessity for precise numerical handling of wave phenomena at each time step. In our numerical experiments, we concentrate on solving the following equations: the constant coefficient Poisson equation, the variable coefficient Poisson equation, the diffusion equation, and the Helmholtz equation. These equations encompass both temporal and spatial derivatives, resulting in a sequence of time-dependent linear systems. At each time step, a sparse linear system akin to that of an elliptic PDE must be tackled. Methods for solving such linear systems are broadly classified into direct and iterative approaches. Direct methods are favored for their high accuracy, stability, and predictability. However, their applicability is constrained when dealing with large-scale systems due to the prohibitive computational and storage requirements. Consequently, iterative methods have emerged as a vital tool for addressing large-scale sparse linear systems. Among these, the algebraic multigrid method (AMG) shines as an effective iterative technique for solving such systems.


In 1987, Ruge, J.W. and St\"{u}ben, K. \cite{ref1} introduced the  { seminal classical} algebraic multigrid (AMG) method, an iterative algorithm tailored for efficiently solving linear algebraic equation systems. This method evolved from the theoretical foundations and principles of general multigrid methods. Unlike traditional methods, AMG does not require knowledge of the geometric and physical properties of the problem; it solely utilizes the information within the coefficient matrix of the linear system to construct a hierarchy of virtual grids. This capability allows multigrid techniques to be applied to a wide variety of problems, including those without a geometric context. The versatility and effectiveness of the AMG method have sparked extensive research interest, encompassing various aspects such as coarsening strategies, smoothers, and practical applications. The AMG parameters used in this paper will be explained in the section on numerical examples.

With the rise of machine learning, researchers have embarked on exploring the integration of AMG with machine learning techniques \cite{ref7}, or the utilization of machine learning to optimize the AMG method \cite{ref9, ref12, ref8}. Our research extends this premise by striving to integrate machine learning methodologies with the AMG approach to augment its performance and applicability.

The effectiveness of the AMG method heavily relies on the selection of smoothing parameters. Traditionally, two main strategies have been employed to determine these parameters: systematic exploration of the parameter space through experimentation or empirical optimization, and theoretical estimation. While the first approach can yield highly accurate optimal parameters, it is time-consuming and often impractical in scientific computing due to the desire for a single, efficient computation. The second strategy, theoretical estimation, may provide a formula or algorithm for computing parameters, but it is problem-specific and any inaccuracies can significantly reduce algorithmic efficiency. Furthermore, as the size of the linear system increases, the applicability of theoretical methods may decrease, especially when it comes to selecting smoothing parameters. To address these challenges, we adopt the approach outlined in \cite{ref16} and utilize Gaussian process regression (GPR) \cite{ref30} to predict optimal parameters. We achieve this by training on a dataset derived from a smaller-scale system, which allows us to discern the mapping from the dimensions of the linear system to more optimal splitting parameters. This strategy avoids the costly parameter traversal and offers an efficient way to forecast relatively optimal parameters within practical computational settings. For the kernel function selection in the GPR method, we incorporate the kernel learning technique described in \cite{ref17, C-2006}. This involves constructing a library of basic kernel functions and formulating the necessary kernel function as a linear combination of these components. The coefficients of this combination are determined through {  Total Least Squares(TLS) training}. This approach not only simplifies the kernel function selection process but also enhances the accuracy of the predictive outcomes.

We have noticed some methods for evaluation metrics \cite{ref24, ref25, ref26, ref27, ref28, ref29}, but they are merely based on one or a few evaluation indicators. Building upon this foundation, our research endeavors to address the limitations observed in prior studies regarding the evaluation of parameter prediction.We transcend the reliance on conventional metrics such as Mean Squared Error (MSE), Root Mean Squared Error (RMSE), and Mean Absolute Error (MAE) by introducing a suite of additional statistical measures. Specifically, we incorporate {  the coefficient of determination} $R^2$ score, also known as the coefficient of determination, which assesses the goodness of fit of our regression model and indicates the proportion of the dependent variable that is predictable from the independent variables. Furthermore, we utilize the correlation coefficient, a statistical measure that quantifies the strength and direction of the linear relationship between two variables.Additionally, we employ the Leave-One-Out cross-validation prediction error (LOO-SPE) to gauge the uncertainty associated with our model predictions. This method provides a robust assessment of model performance by evaluating predictions on each observation individually while using the remaining observations to train the model. Lastly, we incorporate the Bayesian Information Criterion (BIC) value as a scoring criterion, which is useful for model selection by penalizing model complexity and preventing overfitting. By utilizing these comprehensive evaluation metrics, we are well-equipped to thoroughly assess the impact of parameter prediction across various combinations of kernel functions. This approach not only enhances the precision of our parameter predictions but also bolsters their reliability. 
{ Our research contributes to the field by introducing the first kernel learning-based approach for predicting the parameter $\theta$ in AMG, thereby providing a more nuanced understanding of parameter prediction and paving the way for improved algorithmic performance and efficiency.}

{ This paper proposes a GPR-based approach for 
predicting and optimizing parameters in algebraic multigrid 
methods.The effectiveness of the proposed parameter selection 
strategy is demonstrated through numerical experiments on a 
series of PDEs. The paper is structured as follows: Section 2 
outlines the model and algebraic multigrid method and details 
the GPR methodology for parameter prediction, Section 3 
presents the 
results of numerical experiments on various PDEs, and Section 
4 provides a comprehensive summary of our findings.}

\section{Preliminary }\label{sec3}

The performance of the AMG method is closely related to the strong threshold parameter $\theta$, which influences the final computational results and computational cost by determining the C/F splitting. Identifying appropriate parameters for this purpose poses a significant challenge. In this section, we first introduce the basic components of the AMG method, and then we present a novel approach to conduct parameter selection using the GPR method. This method offers a data-driven alternative designed to reduce the computational overhead and enhance the reliability of parameter selection for the AMG method.

\subsection{AMG algorithm}
AMG methods are among the most efficient solvers for linear systems, whose core lies in alternately iterating between coarse and fine grids to handle different components of the error. Taking a two-grid variant as an example: the process begins with performing a predefined number of pre-smoothing iterations on the fine grid to reduce high-frequency error components, followed by transferring the residual from the fine grid to the coarse grid via an restriction operator, and interpolating an error correction, cheaply computed on the coarse grid, back to the fine grid through a prolongation operator to correct the iterative vector, which is finally updated via applying a number of post-smoothing operations. This approach significantly accelerates convergence compared to single-grid methods.

The process can be divided into two stages: SETUP  and SOLVE in AMG. In the SETUP phase, first input the strong threshold parameter $\theta$, which is used to define the strong dependency set of each variable as follows:

\begin{definition} { { Let $N_{i}=\{j|a_{ij}\neq 0 \}$ denote the set of all points that have a non-zero connection to variable  $i$.}}
Given a strong threshold parameter $0<\theta \leq 1$, the variable $i$ strongly depends on the variable $j$ if
\begin{align*}
 |a_{ij}|\geq \theta   \underset{k\neq i}{\rm max}|a_{ik}|,
\end{align*}
{ where $a_{ij}$ is the element in the $i$-th row and $j$-th column of the coefficient matrix $A$.
}
Then the strong dependency set $S_{i}$ and strong influence set $S^{T}_{i}$ of variable $i$ are defined by
\begin{align*}
  S_{i}=\{ j\neq i:~|a_{ij}|\geq { \theta} \underset{k\neq i}{\rm max}|a_{ik}| \} ~~\mbox{and}~~ S^{T}_{i}=\{j|~i\in S_{j}\}.
\end{align*}
\end{definition}

Using criteria $S_{i}$ and $S^{T}_{i}$ for $i=1,2,\cdots,n$ and a specific coarsening strategy , we partition the variable set into a coarse set $C$ and a fine set $F$. { {Common coarsening strategies include the  Ruge-St\"uben algorithm, Cleary–Luby–Jones–Plassman (CLJP) algorithm, and the Parallel Modified Independent Set (PMIS) algorithm . For all AMG experiments presented in this paper, the PMIS algorithm is exclusively adopted as the coarsening strategy.  Further, after coarsening, the prolongation operator $P$ and the restriction operator  $R$ are constructed through a specific interpolation technique, enabling the transfer of residuals and corrections across the grid hierarchy.  Common interpolation techniques include: classical interpolation, standard interpolation, and extended+i (ext+i) interpolation. In the AMG experiments conducted in this paper, the interpolation technique employed is exclusively  ext+i. The specific formula for calculating its interpolation weights is as follows:

\begin{align*}
w_{ij} = -\frac{1}{\tilde{a}_{ii}}
\left(
a_{ij} +
\sum_{k \in F_i^s}
a_{ik}\,
\frac{\bar{a}_{kj}}{\sum_{l \in \hat{C}_i \cup \{i\}} \bar{a}_{kl}}
\right),
\quad j \in \hat{C}_i,
\end{align*}

where
\begin{align*}
\tilde{a}_{ii} = a_{ii}
+ \sum_{n \in N_{i}^w \setminus \hat{C_{i}}} a_{in}
+ \sum_{k \in F_i^s}
a_{ik}\,
\frac{\bar{a}_{ki}}{\sum_{l \in \hat{C_{i}} \cup \{i\}} \bar{a}_{kl}},
\end{align*}
and  $\hat C_{i}$ denotes the extended set of coarse points used for interpolation, and $F_{i}^{S}/C_{i}^{S}$  represents the set of $F/C$-points that strongly   influence variable $i$, $N_{i}^{\omega}=N_{i}\backslash (C_{i}^{S}\cap F_{i}^{S})$    . For details, refer to \cite{D-2008}.}}

During the SOLVE phase, we cycle recursively in V-, W-, F-, K- or any other scheme between different grids the components obtained from the SETUP phase. First, the pre-smoothing process is performed on the fine grid; subsequently, the obtained residual is transferred to the coarse grid via the restriction operator for an exact solution on the coarsest grid or an approximate recursive solution from the hierarchy of coarser grids; then, the corrected error solution is transferred back to the fine grid via the prolongation operator; finally, the resulting solution vector is updated via the post-smoothing process. {  { The AMG component that performs the smoothing operation during either the pre-smoothing process or post-smoothing process is referred to as the smoother. Smoothers are generally classical one-level iterative algorithms for solving linear systems, such as the Jacobi method, Gauss-Seidel method, Symmetric Gauss-seidel method, and Chebyshev method. Throughout this study, the Gauss-Seidel method serves as the smoother in all AMG experiments.}}

Finally,    We present the above two stages in Algorithm \ref{stup} and Algorithm \ref{solve}, respectively.
\begin{small}
	\begin{algorithm}
	\caption{SETUP phase of a two-grid variant}
	\label{stup}
	\begin{enumerate}
		\item \textbf{ Coarsen:} Let $\Omega$ be the set containing all fine-level variables.
		Based on the input $\theta$, construct $S_{i}$ and $S^{T}_{i}$ for each variable $i$ and partition the set $\Omega$ into coarse variables $C$ and fine variables $F$.
		\item \textbf{ Compute $A_{c}$:} Based on the coarse variable set $C$, compute the { prolongation} operator $P$ and restriction operator $R$. Then compute the coarse-level matrix $A_{c}$ by $A_{c}=RAP$.
	\end{enumerate}
\end{algorithm}
\end{small}

\begin{small}
	\begin{algorithm}
	\caption{SOLVE phase of a two-grid variant}
	\label{solve}
	\begin{enumerate}
		\item \textbf{ Pre-smoothing:} smooth $Ax=b$, $k_{1}$ iterations, to get  $x^{k_{1}}$.
				
		\item \textbf{ Coarse-grid correction:}

		{Restrict residual into coarse grid}: $b_{c}=R(b-Ax^{k_1})$\\
		{Obtain $e^{\star}$ by solving the coarse-grid residual equation} $A_{c}e=b_{c}$\\
		{Interpolate and correct}: $x^{\rm new}=x^{k_{1}}+Pe^{*}$
	
			\item \textbf{ Post-smoothing:}  Perform $k_{2}$ iterations on the equation $Ax=b$ using $x^{\rm new}$ as the initial guess
		
	\end{enumerate}
\end{algorithm}
\end{small}


\newpage

\subsection{Gaussian Process Regression}

The training dataset is denoted as 
$ D = \left\{ \left( n_i, \theta_i \right) \mid i = 1, 2, \ldots, d \right\}:= \{ {\boldsymbol{n }}, {\boldsymbol{\theta }} \} $, 
where each pair \((n_i, \theta_i)\) is an input-output correspondence. In this context, \( n_i \) denotes the number of partitions, and \( \theta_i \) signifies the connectivity parameters. Under the assumption that \( \theta_i \) follows a Gaussian Process (GP) relative to \( n_i \), the function evaluations \( f(\boldsymbol{n}) = [f(n_1), f(n_2), \ldots, f(n_d)] \) are collectively distributed as a \( d \)-dimensional joint Gaussian distribution (GD).

\begin{center}
	$
	\left[ \begin{array}{c}
		f\left( n_1 \right)\\
		\vdots\\
		f\left( n_d \right)\\
	\end{array} \right] \thicksim N\left( \left[ \begin{array}{c}
		\mu \left( n_1 \right)\\
		\vdots\\
		\mu \left( n_d \right)\\
	\end{array} \right] ,\left[ \begin{matrix}
		k\left( n_1,n_1 \right)&		\cdots&		k\left( n_1,n_d \right)\\
		\vdots&		\ddots&		\vdots\\
		k\left( n_d,n_1 \right)&		\cdots&		k\left( n_d,n_d \right)\\
	\end{matrix} \right] \right) .
	$
\end{center}
The Gaussian process is defined by its mean function \(\mu(n)\) and covariance function \(k(n, n')\). We represent the GP in the following manner:
\[
f(n) \sim GP(\mu(n), K(n, n')),
\]
where \(n\) and \(n'\) represent any two elements within the input set \(\boldsymbol{n}\). Conventionally, for the sake of simplicity, we set the mean function to zero.

The goal of Gaussian process regression is to establish a mapping between the input set {$\boldsymbol{n}$} and the output set {$\boldsymbol{\theta }$}, represented as \( f(n): \mathbb{R} \rightarrow \mathbb{R} \), and to forecast the potential output value \(\theta_* = f(n_*)\) at a new test point \(n_*\). In practical linear regression scenarios, our model is given by:
\[
\theta = f(n) + \eta,
\]
where \(\theta\) is the observed value that incorporates noise \(\eta\). 
{  We presume that \(\eta\) follows a Gaussian distribution with a zero mean and variance \(\sigma^2\), that is \(\eta \sim N(0, \sigma^2)\).}
 The common range for \(\sigma\) is \([10^{-6}, 10^{-2}]\). For this investigation, we have set \(\sigma = 10^{-4}\).

Then, we can express the prior distribution of the observed values $\boldsymbol{\theta}$ as
\begin{center}
	$
	\boldsymbol{\theta }\thicksim N\left( \boldsymbol{\mu }_{\theta}\left( \boldsymbol{n} \right) ,K\left( \boldsymbol{n},\boldsymbol{n} \right) +\sigma ^2I_d \right).
	$
\end{center}
The joint prior distribution of the observed value $\theta$ and the predicted value $\theta_*$ is given by
\begin{center}
	$
	\left[ \begin{array}{c}
		\boldsymbol{\theta }\\
		\boldsymbol{\theta }_*\\
	\end{array} \right] \thicksim N\left( \left[ \begin{array}{c}
		\boldsymbol{\mu }_{\theta}\\
		\boldsymbol{\mu }_{\theta _*}\\
	\end{array} \right] ,\left[ \begin{matrix}
		K\left( \boldsymbol{n},\boldsymbol{n} \right) +\sigma ^2I_d&		K\left( \boldsymbol{n},\boldsymbol{n}_* \right)\\
		K\left( \boldsymbol{n}_*,\boldsymbol{n} \right)&		K\left( \boldsymbol{n}_*,\boldsymbol{n}_* \right)\\
	\end{matrix} \right] \right) ,
	$
\end{center}
where $I_d$ is the $d$-dimensional identity matrix, $K(\boldsymbol{n},\boldsymbol{n})=(k_{ij})$ is the symmetric positive definite covariance matrix, and $k_{ij}=k(n_i,n_j)$. $K(\boldsymbol{n},\boldsymbol{n}_*)$ is the symmetric covariance matrix between the training set $\boldsymbol{N}$ and the test set $\boldsymbol{n}_*$.

According to Bayes formula, the posterior distribution of the predicted value \(\boldsymbol{\theta}_*\) given the observed values \(\boldsymbol{\theta}\) is:
\begin{equation}
	p\left( \boldsymbol{\theta }_*|\boldsymbol{\theta } \right) =\frac{p\left( \boldsymbol{\theta }|\boldsymbol{\theta }_* \right) p\left( \boldsymbol{\theta }_* \right)}{p\left( \boldsymbol{\theta } \right)},
	\label{beiyesi}
\end{equation}
where \( p(\boldsymbol{\theta}_*) \) is the prior distribution of the predicted value and \( p(\boldsymbol{\theta}) \) is the marginal likelihood of the observed values.

The joint posterior distribution of the prediction is given by:
\begin{equation}\label{houyan}
	\boldsymbol{\theta }_*|\boldsymbol{n},\boldsymbol{\theta },\boldsymbol{n}_*\thicksim N\left( \boldsymbol{\mu }_*,\boldsymbol{\sigma }_{*}^{2} \right) ,
\end{equation}
where the mean \(\boldsymbol{\mu}_*\) and variance \(\boldsymbol{\sigma}_*^2\) are:
\begin{center}
	$
	\begin{array}{l}
		\boldsymbol{\mu }_*=K\left( \boldsymbol{n}_*,\boldsymbol{n} \right) \left[ K\left( \boldsymbol{n},\boldsymbol{n} \right) +\sigma ^2I_d \right] ^{-1}\left( \boldsymbol{\theta }-\boldsymbol{\mu }_{\theta} \right) +\boldsymbol{\mu }_{\theta _*},\\
		\boldsymbol{\sigma }_{*}^{2}=K\left( \boldsymbol{n}_*,\boldsymbol{n}_* \right) -K\left( \boldsymbol{n}_*,\boldsymbol{n} \right) \left[ K\left( \boldsymbol{n},\boldsymbol{n} \right) +\sigma ^2I_d \right] ^{-1}K\left( \boldsymbol{n},\boldsymbol{n}_* \right) .\\
	\end{array}
	$
\end{center}
For the output \(\theta_*\) in the test set, we use the mean of the Gaussian regression as its estimate, i.e., \(\widehat{\theta}_* = \mu_*\).

The kernel function plays a pivotal role in GPR, as it is responsible for generating the covariance matrix that quantifies the similarity between any two input variables. The closer the distance between inputs, the higher the correlation between their corresponding output variables. Hence, it is essential to select or construct kernel functions that meet the specific demands of the problem at hand. Commonly utilized kernel functions include radial basis functions, rational quadratic kernel functions, exponential kernel functions, and periodic kernel functions.

In this article, we have employed a kernel learning approach. We have established a kernel function library \(\mathcal{K} = \left\{ k_{\xi}(x, x') \mid \xi = 1, \cdots, N \right\}\), which encompasses a variety of fundamental kernel functions and their multiplicative combinations. Subsequently, we derive the kernel function by linearly combining elements from this library:
\[
k(x, x') = \sum_{\xi=1}^N c_{\xi} k_{\xi}(x, x').
\]

For the selection of the hyperparameter \(\alpha\) in the kernel function, we obtain it by maximizing the log-likelihood function:
\[
L = \log p(y | X, \alpha) = -\frac{1}{2} y^T \left( K(X, \alpha) + \sigma^2 I_n \right)^{-1} y - \frac{1}{2} \log \det \left( K(X, \alpha) + \sigma^2 I_n \right) - \frac{n}{2} \log 2\pi
\]
{ where $K(X,\alpha)$ is the $n\times n$  kernel matrix computed from the input dataset $X$ given the hyperparameters $\alpha$, with entires $[K(X,\alpha)]_{ij} = k(x_i,x_j;\alpha)$. Here $\alpha$ includes the linear combination coefficients $c_\xi$ and the internal parameters of the base kernels. By maximizing $L$, we perform empirical Bayes estimation of the hyperparameters $\alpha$. This approach balances data fitting and model complexity penalty , which helps prevent overfitting and select an appropriate kernel combination.}

\subsection{Multi-Indicator Driven Evaluation}

In this section, we will present a comprehensive set of evaluation metrics crafted to assess the performance of GPR, specifically in the domain of parameter prediction. These metrics offer diverse insights into the quality of the predictions, with each one reflecting different aspects of predictive performance based on their values. We will explore the details of each metric and explain how their values are indicative of the prediction model's effectiveness.

\subsubsection{Positive Indicators: Value Enhancement and Prediction Optimization}

In evaluating the performance of parameter predictions in Gaussian process regression, we employ a variety of metrics to thoroughly assess the model's effectiveness. This set of metrics includes the coefficient of determination \( R^2 \), the correlation coefficient (Corr), and the Prediction Interval Coverage Probability (PICP), which together form a robust framework for measuring the model's predictive power.

The coefficient of determination, represented as \( R^2 \), is a statistical measure that indicates the proportion of the variance in the dependent variable that can be explained by the independent variables within the regression model. It offers a measure of how well the observed outcomes are replicated by the model, relative to the total variation present in the data.

The formula for calculating \( R^2 \) is as follows:
\[
R^2 = 1 - \frac{\sum_{i=1}^{n}(y_i - \hat{y}_i)^2}{\sum_{i=1}^{n}(y_i - \bar{y})^2},
\]
where
\( y_i \) is the observed value of the \( i \)-th data point, \( \hat{y}_i \) is the predicted value of the \( i \)-th data point from the model, \( \bar{y} \) is the mean of the observed values.{  The \( R^2 \) value ranges between 0 and 1, where a value of 1 indicates that the model perfectly predicts the data, and a value of 0 suggests that the model does not explain any of the variability in the response data around its mean.}



Corr is a statistical tool that measures the degree and direction of the linear relationship between a model's predicted values and the actual data points\cite{ref18}. 
{  This correlation coefficient is derived by dividing the covariance of the predicted and observed values by the product of their respective standard deviations, resulting in a value that falls within the range of -1 to 1. In this range, values nearing 1 or -1 denote a robust linear association. However, it is crucial to recognize that a high correlation does not imply causation; it signifies only a strong association between the variables. Moreover, the correlation is sensitive to outliers, making it essential to scrutinize the data for extreme values that could disproportionately affect it.}

\[
Corr = \frac{\sum_{i=1}^{n}(x_i - \bar{x})(y_i - \bar{y})}{\sqrt{\sum_{i=1}^{n}(x_i - \bar{x})^2 \sum_{i=1}^{n}(y_i - \bar{y})^2}},
\]
{ where \(x_i\) and \(y_i\) represent the individual sample points indexed by \(i\), \(\bar{x}\) and \(\bar{y}\) denote the means of the \(x\) and \(y\) values, respectively, and \(n\) is the total number of observations. }


GPR not only predicts parameters but also furnishes a quantitative assessment of the uncertainty surrounding these predictions, encapsulated by the PICP \cite{ref19,ref20}. The PICP at a 0.95 confidence level reveals the model's capacity to ensure that true values fall within the predicted intervals. By determining the mean  \(\mu_*(x)\) and variance  \(\sigma_*^2(x)\) via GPR, the 0.95 confidence interval can be computed as:

\[
I_{\theta}^{0.95}(x) = \left[\mu_*(x) - 1.96 \sigma_*(x), \mu_*(x) + 1.96 \sigma_*(x)\right].
\]

Given a test set:

\[
\left\{x_{i}^{\text{test}}, i = 1, \ldots, N\right\} \in \mathbb{X}^N,
\]
where the values of \(z\) are known (actual values), the prediction interval coverage rate is:

\[
PICP = \frac{1}{N} \sum_{i=1}^{N} \mathbf{1}_{z(x_i^{\text{test}}) \in I_{\theta}^{0.95}(x_i^{\text{test}})}.
\]

Here, \(\mathbf{1}_{z(x_i^{\text{test}}) \in I_{\theta}^{0.95}(x_i^{\text{test}})}\) is an indicator function that equals 1 if the actual value \(z(x_i^{\text{test}})\) falls within the 95\% confidence interval \(I_{\theta}^{0.95}(x_i^{\text{test}})\), and 0 otherwise.

The {PICP} provides a measure of how frequently the model's predicted intervals contain the true values. A high PICP value indicates that the model's predicted intervals are reliable and include the true values a significant proportion of the time. This is a useful metric for evaluating the reliability of the model's predictions and for establishing appropriate confidence levels for the predictions.

{ 
The combined use of these performance metrics provides a comprehensive evaluation of the predictive prowess of the GPR model. High \( R^2 \) and Corr indicate that the model effectively captures the underlying trends in the data. Simultaneously, a strong PICP underscores the reliability and consistency of the model's predictive outputs.} 

\subsubsection{Negative Indicators: Value Reduction and Error Decrease}

In the domain of assessing parameter predictions from GPR models, a suite of negative indicators is utilized to measure the magnitude of prediction errors. This collection includes Mean Squared Error (MSE), Root Mean Squared Error (RMSE)\cite{ref21}, Mean Absolute Error (MAE), Bayesian Information Criterion (BIC), Median Absolute Percentage Error (MdAPE), and Leave-One-Out Cross-Validated Score Error (LOO-SPE). Collectively, these metrics form a comprehensive framework for evaluating the model's predictive accuracy.

{ MSE, RMSE, and MAE} serve as direct measures of prediction error, reflecting the discrepancies between predicted values and actual observations from various perspectives. MSE assesses the model's precision by averaging the squared differences between predictions and actual values, providing a clear view of the average squared error. RMSE offers a scaled measure of error magnitude, expressed in the same units as the data, which is particularly useful for grasping the practical significance of errors. MAE provides a robust measure of error by averaging the absolute differences between predictions and actual values, making it less susceptible to the influence of outliers. In an ideal scenario, these metrics should yield lower values, signifying a higher level of predictive precision for the model.

The Bayesian Information Criterion (BIC) is a pivotal tool in model selection that imposes a penalty on model complexity to mitigate the risk of overfitting. The BIC is formulated to consider both the maximum likelihood estimate of the model and the number of parameters it contains; a lower BIC value indicates a more favorable balance between the model's fit to the data and its simplicity.The formula for calculating the BIC is given below:

\[
\mathrm{BIC} = -2 \ln (L) + k \ln (n),
\]
where \(-2 \ln(L)\)  incentivizes the model to achieve a high likelihood value, signifying a good fit to the data. Conversely, the term \(k \ln(n)\) increases with the number of parameters \(k\) in the model, thus imposing a penalty for increased complexity. This aspect of the BIC ensures that models with an excessive number of parameters are penalized, even if they fit the data well, as their BIC values will be elevated. Consequently, when evaluating different models, the one exhibiting a smaller BIC value is often deemed more desirable, as it strikes a better trade-off between data fit and model simplicity.


The Median Absolute Percentage Error (MdAPE) is a robust metric for assessing the predictive accuracy of a forecasting model\cite{ref22}. It measures the median of the absolute percentage differences between the predicted and actual values. The MdAPE is calculated as follows:

\[
\operatorname{MdAPE} = \operatorname{median}\left(\left|\frac{x_i - \hat{x}_i}{x_i}\right|\right),
\]
where, \(x_i\) denotes the true value of the \(i\)-th sample, and \(\hat{x}_i\) denotes the predicted value of the \(i\)-th sample. The Median Absolute Percentage Error (MdAPE) utilizes the median instead of the mean, which makes it less sensitive to outliers, thereby offering a more dependable measure of the model's predictive performance in scenarios where extreme values are present.


{ The Median Absolute Percentage Error (MdAPE) is a key metric for evaluating predictive model accuracy, where a lower value indicates higher precision. Conventionally, an MdAPE below 0.1 signifies excellent performance, while a value above 0.5 suggests poor performance. }

In addition to the previously discussed metrics, we introduce an evaluation method based on a scoring criterion known as LOO-SPE \cite{ref23}, which is used to assess the predictive performance of the GPR model. This method is characterized by its focus on minimizing the prediction error through the Leave-One-Out cross-validation approach. The lower the LOO-SPE value, the better the model's predictive performance is considered to be.

Within the framework of the GPR model, the predictive distribution \( P_{\theta,-i} \) is characterized by the following conditional Gaussian distribution:

\[
P_{\theta,-i} { =} \mathcal{N}\left(\mu_{\theta
where,-i}, \sigma_{\theta,-i}^{2}\right),
\]

\[
\mu_{\theta,-i} = y_{i} - \frac{\left(C(X, X)^{-1}(y - m)\right)_{i}}{\left(C(X, X)^{-1}\right)_{i, i}} \quad \text{and} \quad \sigma_{\theta,-i}^{2} = \frac{1}{\left(C(X, X)^{-1}\right)_{i, i}}.
\]

Here, \(\mu_{\theta,-i}\) represents the predictive mean, \(\sigma_{\theta,-i}^{2}\) denotes the predictive variance, and \(\theta\) signifies the model parameters.

In the Leave-One-Out method, for each training sample \( y_{i} \), a conditional model \( P_{\theta,-i}(y \mid x) \) is constructed, trained on the dataset with the \( i \)-th sample excluded. The LOO-SPE scoring criterion is articulated as:

\[
J_{n}^{S}(\theta) = \frac{1}{n} \sum_{i=1}^{n} S\left(P_{\theta,-i}, y_{i}\right),
\]
where \( S \) is the scoring rule, and \( P_{\theta,-i} \) denotes the distribution of \( Y(x_i) \) given the condition \( Y(x_j) = y_j \) for \( j \neq i \).

We define the mean of the predictive distribution \( P \) as \(\mu\), and consider the deviation of \( y \) from \(\mu\). The squared prediction error (SPE) is given by:

\[
S^{\mathrm{SPE}}(P, y) = (y - \mu)^{2}.
\]

Specifically, integrating the results from the GPR model, the LOO-SPE can be articulated as:

\[
J_{n}^{\mathrm{SPE}}(\theta) = \frac{1}{n} \sum_{i=1}^{n} \left(\mu_{\theta,-i} - y_{i}\right)^{2} = \frac{1}{n} \sum_{i=1}^{n} \left(\frac{\left(C(X, X)^{-1}(y - m)\right)_{i}}{\left(C(X, X)^{-1}\right)_{i, i}}\right)^2.
\]

In practical applications, relying on a single metric often falls short of capturing the entire scope of a model's performance. Hence, this paper adopts a comprehensive evaluation strategy that utilizes multiple indicators to achieve a more holistic assessment of predictive performance. In the forthcoming section on numerical experiments, we will detail how these indicators can be synthesized and how they collectively enhance our understanding of the GPR model's performance in parameter prediction.

\section{Numerical experiments}

In this section, we focus on the role of GPR in parameter selection and provide various numerical examples to demonstrate its effectiveness. These examples not only cover a broad range of mathematical models but also showcase the performance and adaptability of GPR in various scenarios through experiments. We hope that readers can intuitively grasp the practical application value and potential advantages of GPR.
To ensure the precision and credibility of the experimental results, all numerical experiments described in this section are performed by discretizing the equations using the finite element method in the IFEM software package \cite{CL2008}. If the number of partitions is $n$, the order of the final discretized matrix is $(n+1)^{2}$.  The discretized systems are then solved on the hypre platform \cite{F-Y-2022}.
Before the experiments begin, we perform a uniform triangulation of the computational domain to construct an appropriate numerical model.
In these experiments, the  smoother is  Gauss-Seidel , with and PMIS as the coarsening algorithm.

We define the termination condition for our iterative process as follows:
\[
tol^k = \frac{\left\| b - Ax^k \right\|_2}{\left\| b \right\|_2}.
\]
Here, \( tol^k \) denotes the tolerance at the \( k \)-th iteration, \( b \) represents the known vector, \( A \) is the non-singular matrix, and \( x^k \) signifies the approximate solution at the \( k \)-th iteration. The initial tolerance is set to \( 0 \). The iterative process is halted when \( tol^k < 1 \times 10^{-08} \) or when the iteration count reaches a pre-established maximum threshold, which is typically set to 200 unless otherwise indicated.
And the column “Speed up” is the CPU time of default $\theta$ divided by the CPU time of predicted $\theta$.



When choosing the best parameters, we focus on decreasing iterations to boost computational efficiency and numerical solution stability. If different sets yield the same iteration count, we opt for the one with the lowest error, ensuring a balance between efficiency and precision. This section uses a V-cycle approach to address linear systems.

Our numerical examples aim to highlight the practicality and performance of our GPR-based parameter selection method across various equations. We aim to showcase the GPR method's advantages in optimization, validate its theoretical basis, and underscore its robustness and effectiveness in facing real-world complexities.

\subsection{Constant Coefficient Poisson Equation}

We consider the constant coefficient Poisson equation, a fundamental partial differential equation in mathematical physics, given by:
\[
\begin{cases}
	-\nabla \cdot \left( a \nabla u \right) = f, & (x,y) \in \Omega = [0,1]^2, \\
	u = 0, & (x,y) \in \partial \Omega.
\end{cases}
\]
Here, the coefficient \( a = 1 \), and the exact solution is \( u = \sin(\pi x)\sin(\pi y) \). The right-hand side \( f \) is derived from the exact solution by substituting \( u \) into the Poisson equation. To find the optimal parameters for our GPR-based parameter selection method, we perform a grid search within the interval [0,1] with a step size of 0.001. The number of partitions \(n \) ranges from 64 to 400 with an increment of 16. {  The traversal identified the optimal parameters, and the results are detailed in Table \ref{tab:possion1}}

\begin{table}[h]
	\centering
	\caption{Traversal results of AMG solving the constant coefficient Poisson equation with different $n$}
	\label{tab:possion1}
	\begin{tabular}{|c|c|c|c|c|c|}
		\hline
		$n$ & $\theta$ & iter & $n$ & $\theta$ & iter  \\ \hline
		64 & 0.256 & 13 & 240 & 0.268 &  15 \\ \hline
		80 & 0.269 & 13 & 256 & 0.284 &15 \\ \hline
		96 & 0.387 & 13 & 272 & 0.267& 15 \\ \hline
		112 & 0.324& 13 & 288 &0.287 &15 \\ \hline
		128 & 0.222& 14 & 304& 0.311& 15 \\ \hline
		144 & 0.246& 14 & 320 &0.275& 15 \\ \hline
		160 & 0.293& 14 & 336 &0.285 & 15 \\ \hline
		176& 0.352& 14 &352 &0.319& 15 \\ \hline
	    192 &0.351 &14& 368 &0.314& 15 \\ \hline
		208 & 0.334 &14 & 384 &0.349& 15 \\ \hline
		224 & 0.232 & 15 & 400& 0.34& 16 \\ \hline
	\end{tabular}
\end{table}

\begin{figure}[H]
	\centering
	\label{fig:possion1}
	\caption{The variation of the number of iterations for AMG solving the constant coefficient Poisson equation with different connectivity parameters $\theta$ when $n$ is 64, 256, and 400, respectively}
	\begin{minipage}{0.32\linewidth}
		\centering
		\includegraphics[width=1\linewidth]{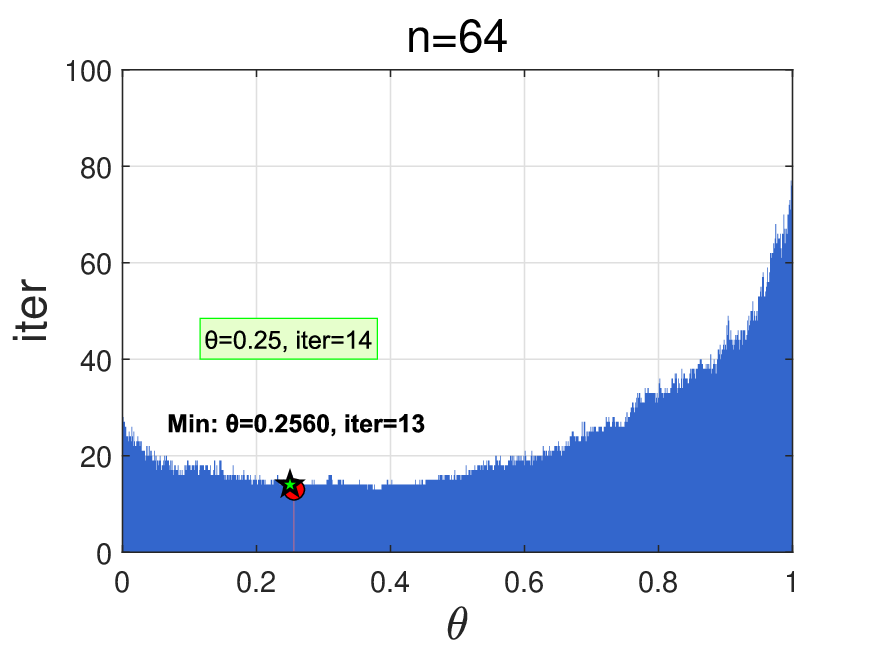}
	\end{minipage}
	\begin{minipage}{0.32\linewidth}
		\centering
		\includegraphics[width=1\linewidth]{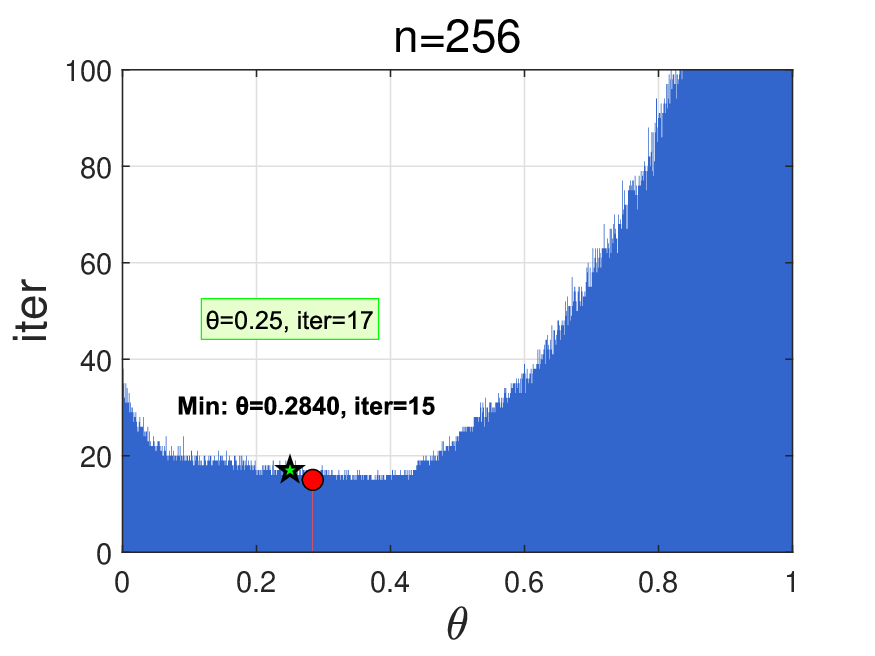}
	\end{minipage}
		\begin{minipage}{0.32\linewidth}
		\centering
		\includegraphics[width=1\linewidth]{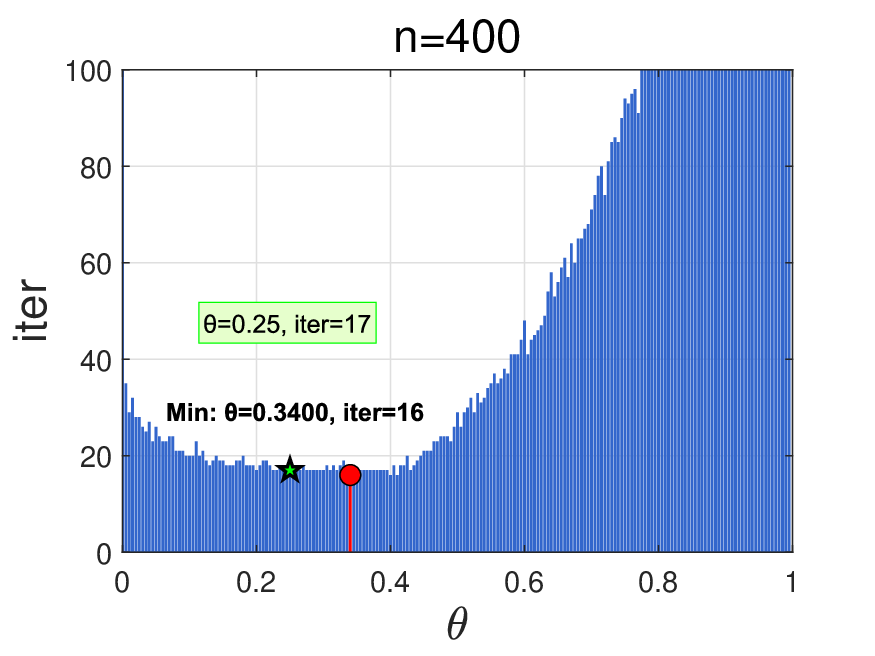}
	\end{minipage}
\end{figure}

Figure \ref{fig:possion1} illustrates the fluctuation in the number of iterations for the AMG method when \(n = {64,~256,~400} \), across various traversal parameters. The traversal for the connectivity parameter \(\theta\) spans the interval \([0,1]\) with a step size of 0.001, and the iteration process is capped at a maximum of 100 iterations. Figure \ref{fig:possion1} clearly demonstrates that the number of iterations fluctuates with changes in the connectivity parameter \(\theta\), and this fluctuation intensifies as \(n \) increases. This observation underscores the sensitivity of the AMG method's efficiency to the choice of the connectivity parameter \(\theta\).

\begin{table}[H]
	\centering
{
	\caption{Training set and retraining set for using GPR to predict the optimal parameters for AMG solving the constant coefficient Poisson equation}
	\label{tab_hmxunlianji}
	\scalebox{1}{
		\begin{tabular}{|c|c|c|}
			\hline
			Training set & $n\in[64, 400], \Delta n=16$\\
			Retraining set 1 & $n\in[200,600]$, 10 randomly selected points\\
			Retraining set 2 & $n\in[200,600]$,12 randomly selected points \\
			\hline
		\end{tabular}
	}
}
\end{table}

We have successfully predicted the optimal parameters for the constant coefficient Poisson equation using GPR method. The following table outlines the composition of the training set, test set, and retraining set used in our GPR model. The training set was constructed by selecting values of \(n \) from 64 to 400, with an increment of \(\Delta n = 16\).
For retraining we adopt targeted sampling: the first retraining phase selects 10 new points uniformly distributed in [200,600], followed by a second retraining phase where 12 additional points are chosen from [200,600] to refine the model further. We summarize this in {  Table \ref{tab_hmxunlianji}}.

\begin{figure}[H]
	\centering
	\caption{Regression curves for predicting $\theta$ with respect to $n$ using GPR for AMG solving the constant coefficient Poisson equation}
	\label{Fig:possion_pre}
	\begin{minipage}{0.32\linewidth}
		\centering
		\includegraphics[width=1\linewidth]{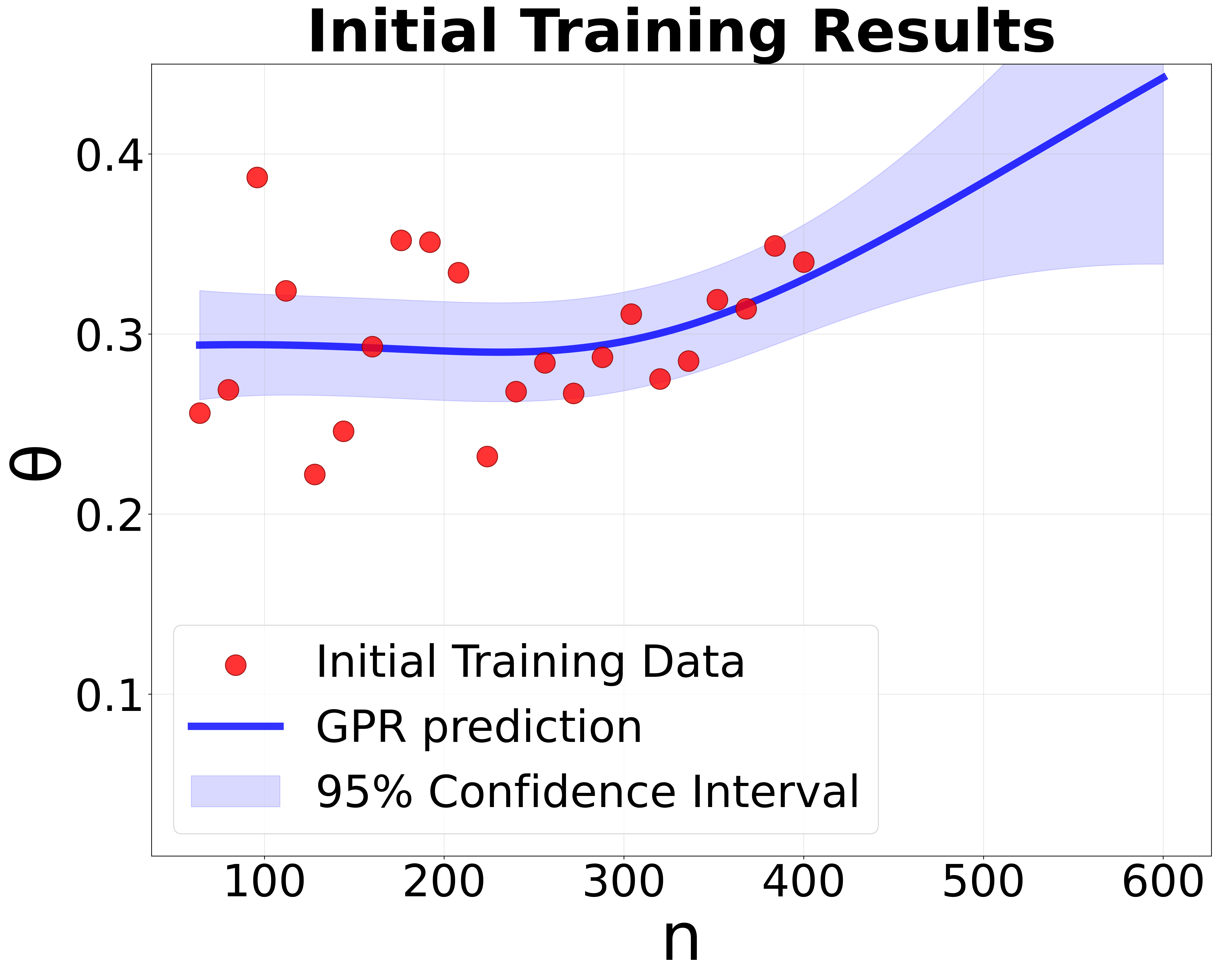}
	\end{minipage}
	\begin{minipage}{0.32\linewidth}
		\centering
		\includegraphics[width=1\linewidth]{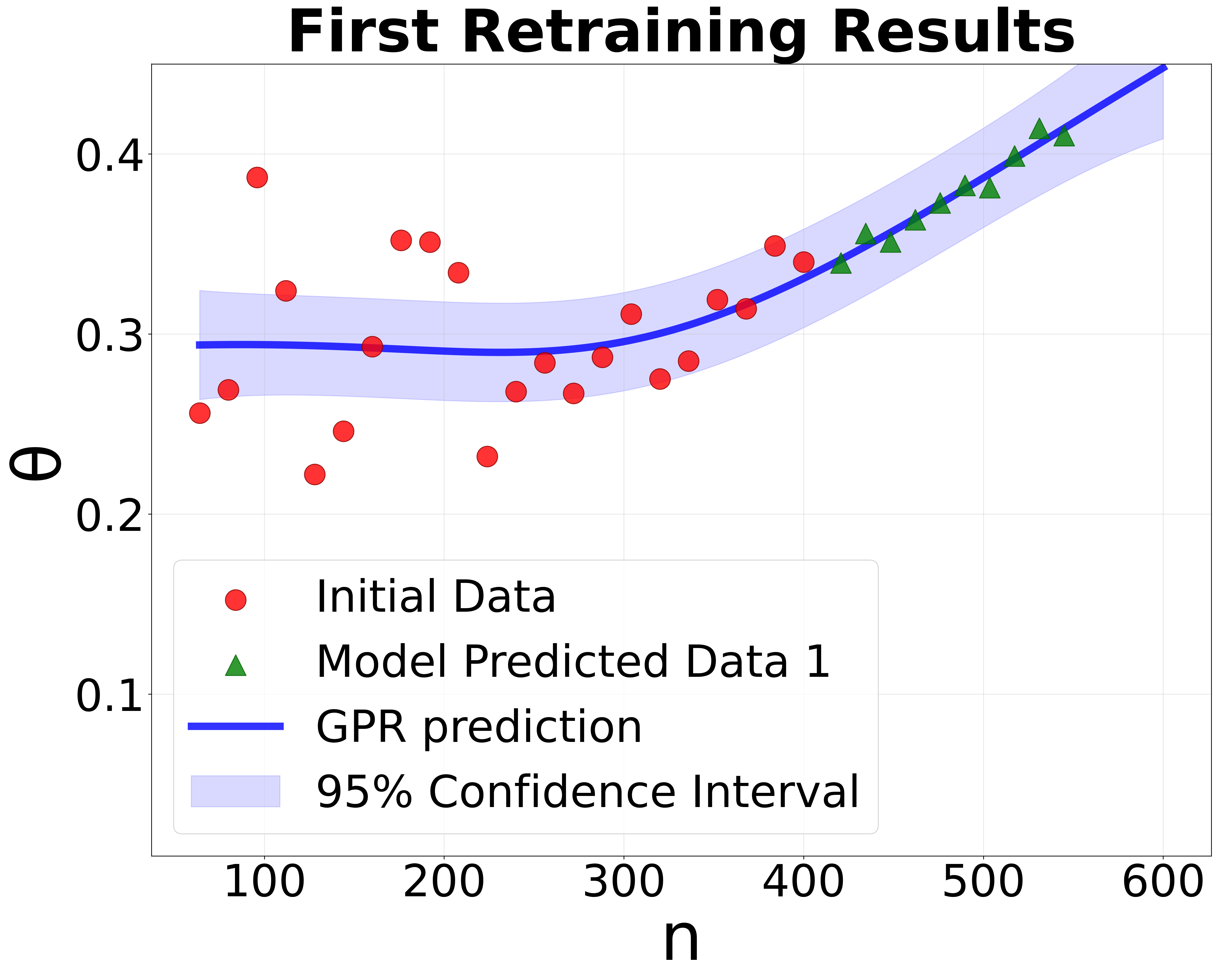}
	\end{minipage}
		\begin{minipage}{0.32\linewidth}
		\centering
		\includegraphics[width=1\linewidth]{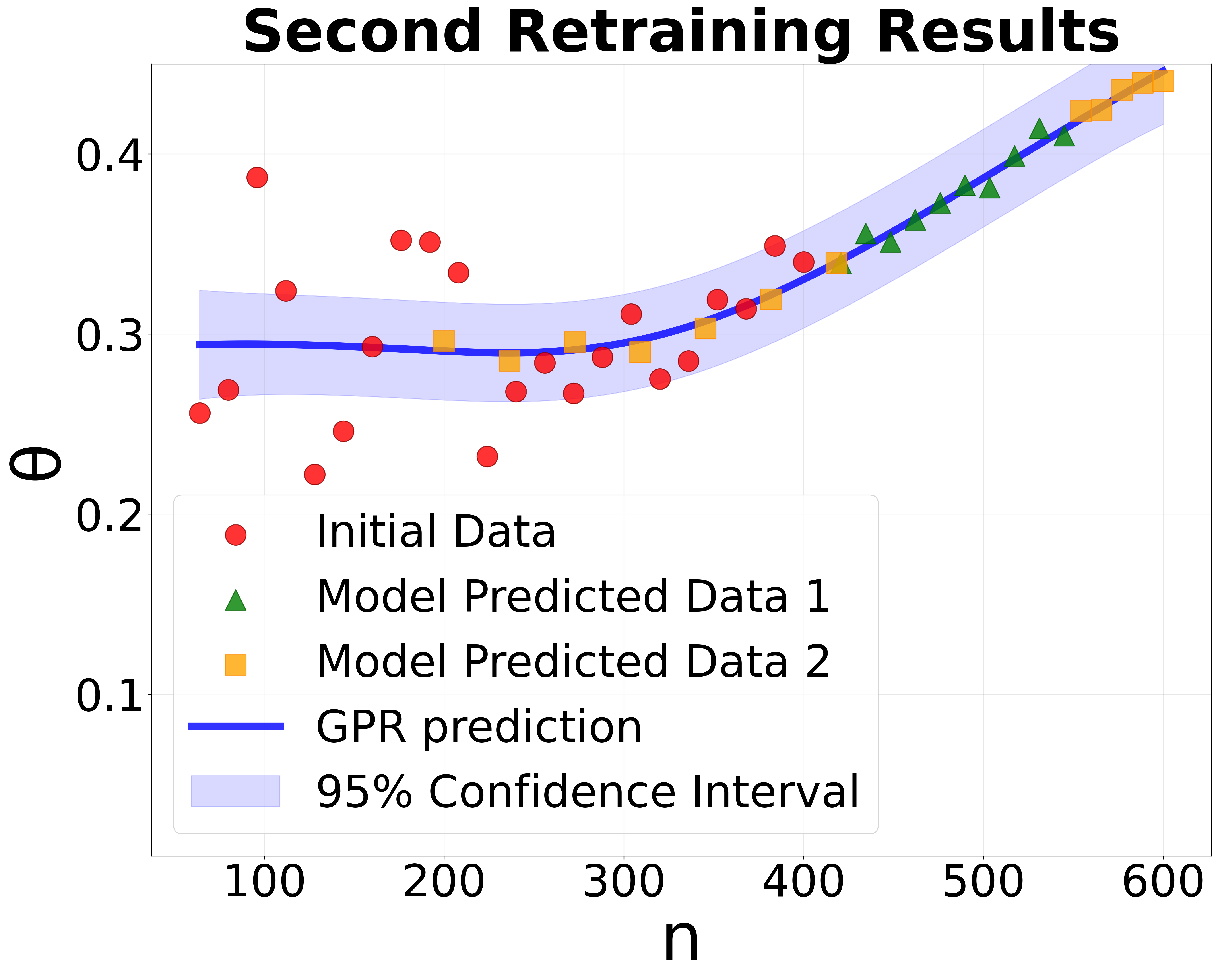}
	\end{minipage}
\end{figure}

\begin{figure}[H]
	\centering
	\caption{{  Variation of the number of iterations for AMG solving the constant coefficient Poisson equation and different connectivity parameters $\theta$ when $n$ is ~ 1024,~ 1200,~1400,~1600,~1800,~2048, respectively}}
	\label{hm_travel_more}
		\begin{minipage}{0.32\linewidth}
		\centering
		\includegraphics[width=1.1\linewidth]{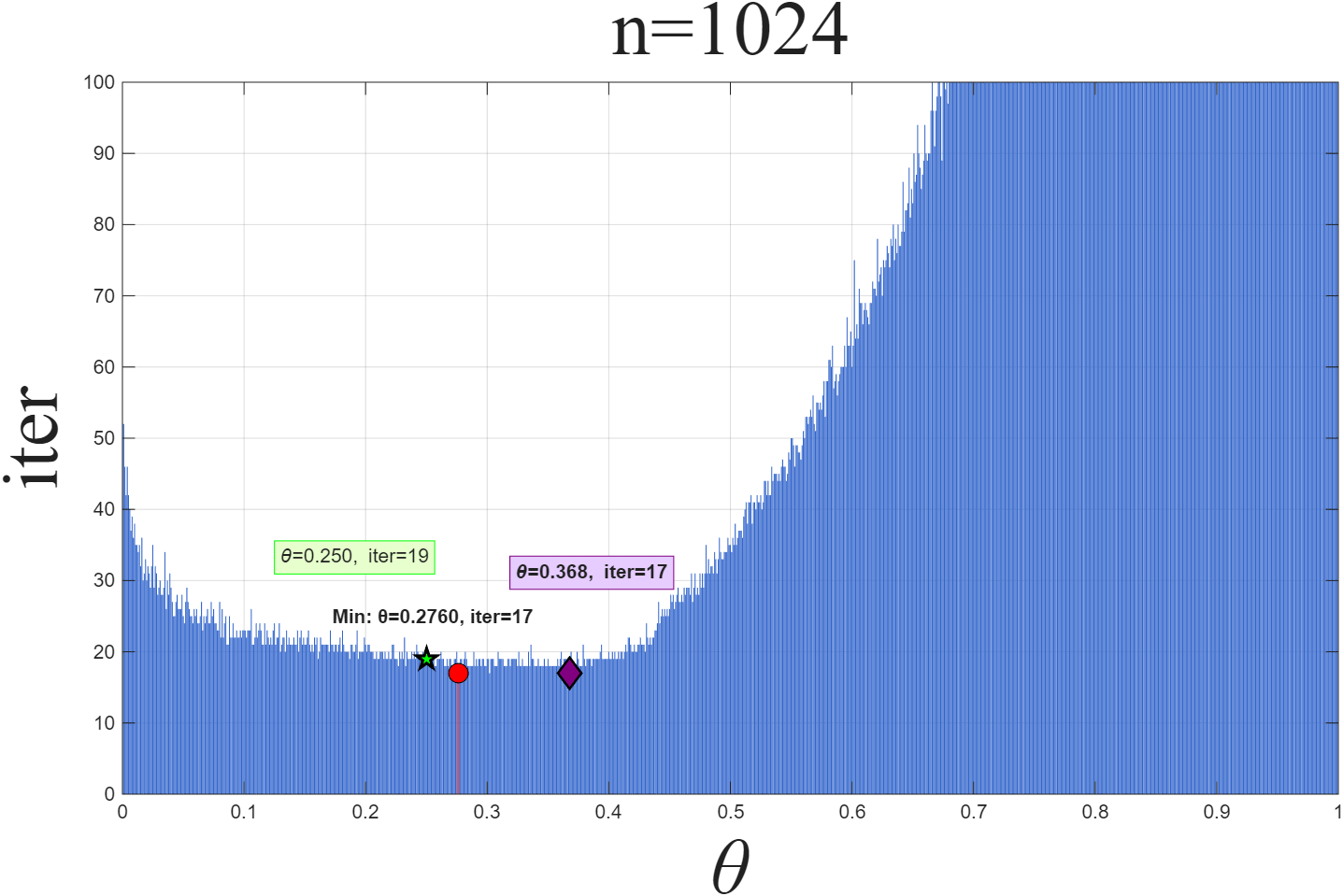}
	\end{minipage}
	\begin{minipage}{0.32\linewidth}
		\centering
		\includegraphics[width=1.1\linewidth]{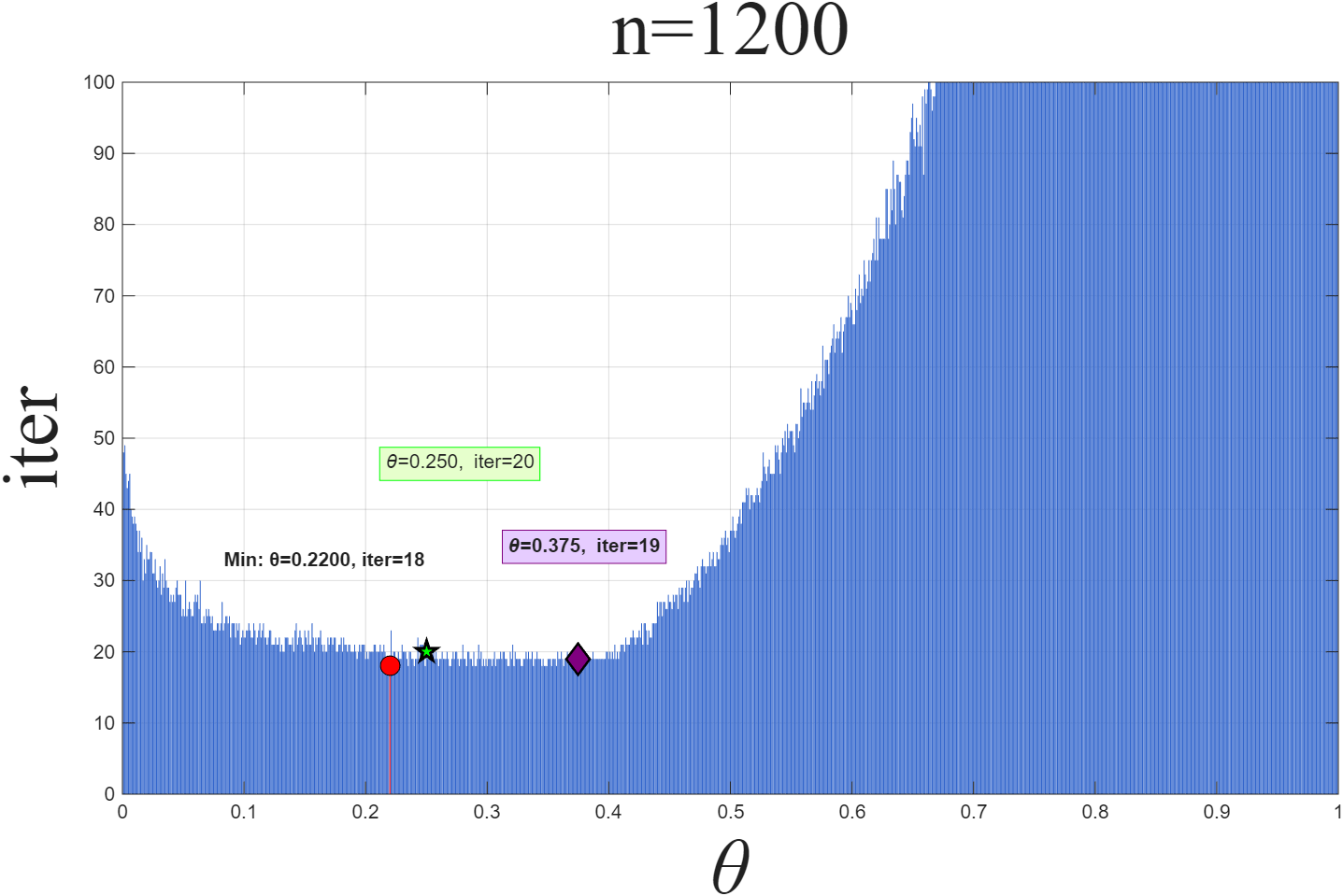}
	\end{minipage}
	\begin{minipage}{0.32\linewidth}
		\centering
		\includegraphics[width=1.1\linewidth]{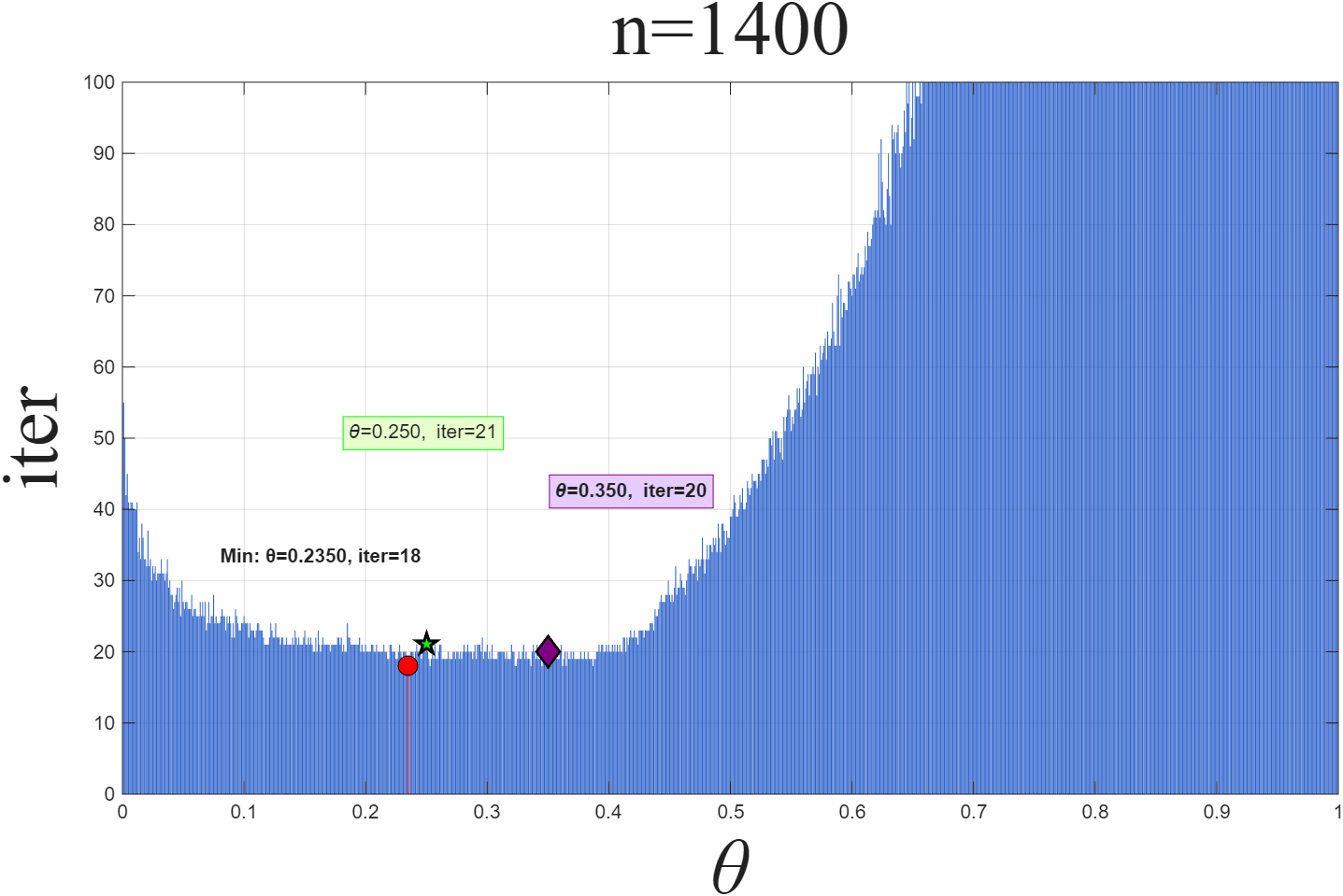}
	\end{minipage}
		\begin{minipage}{0.32\linewidth}
		\centering
		\includegraphics[width=1.1\linewidth]{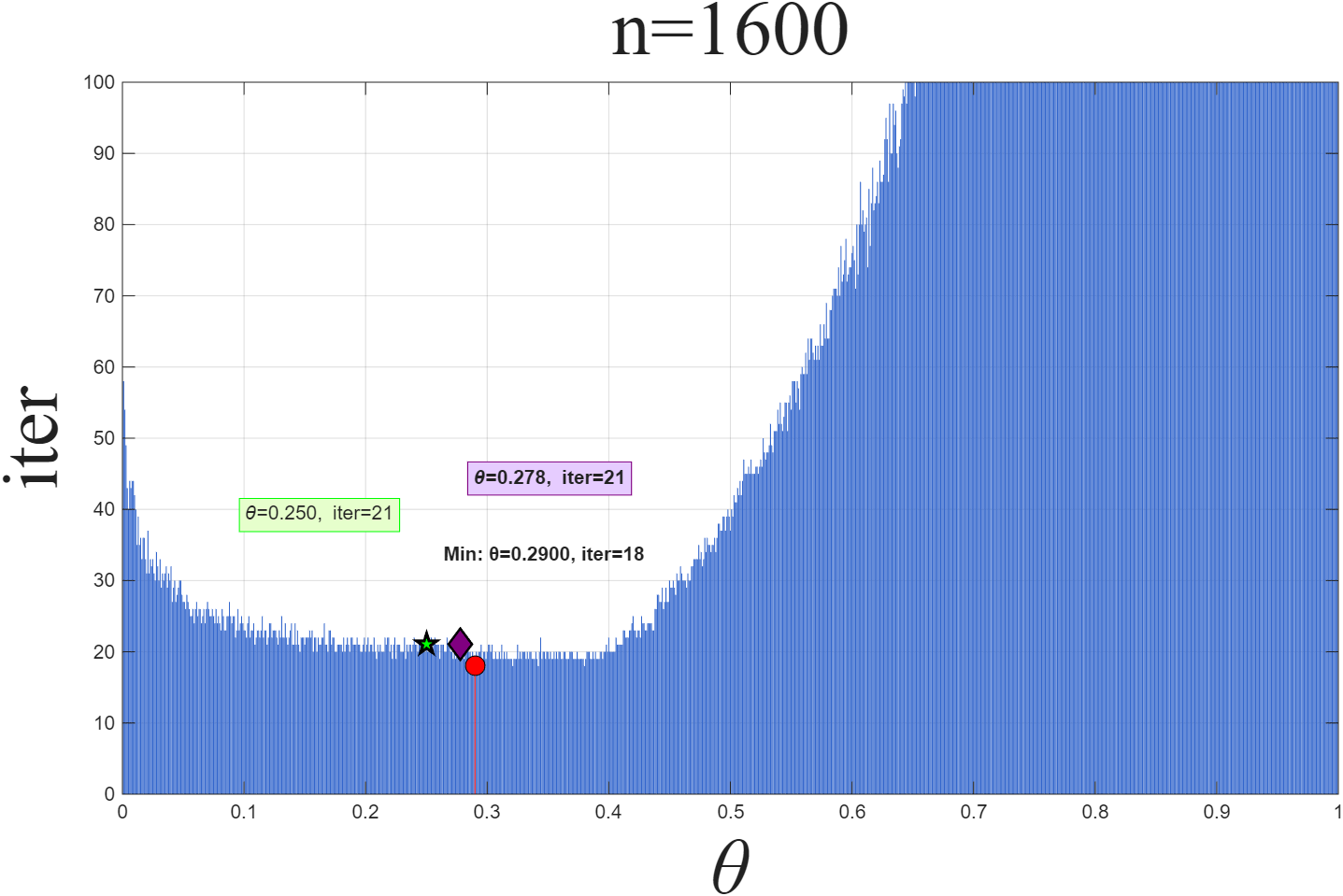}
	\end{minipage}
			\begin{minipage}{0.32\linewidth}
		\centering
		\includegraphics[width=1.1\linewidth]{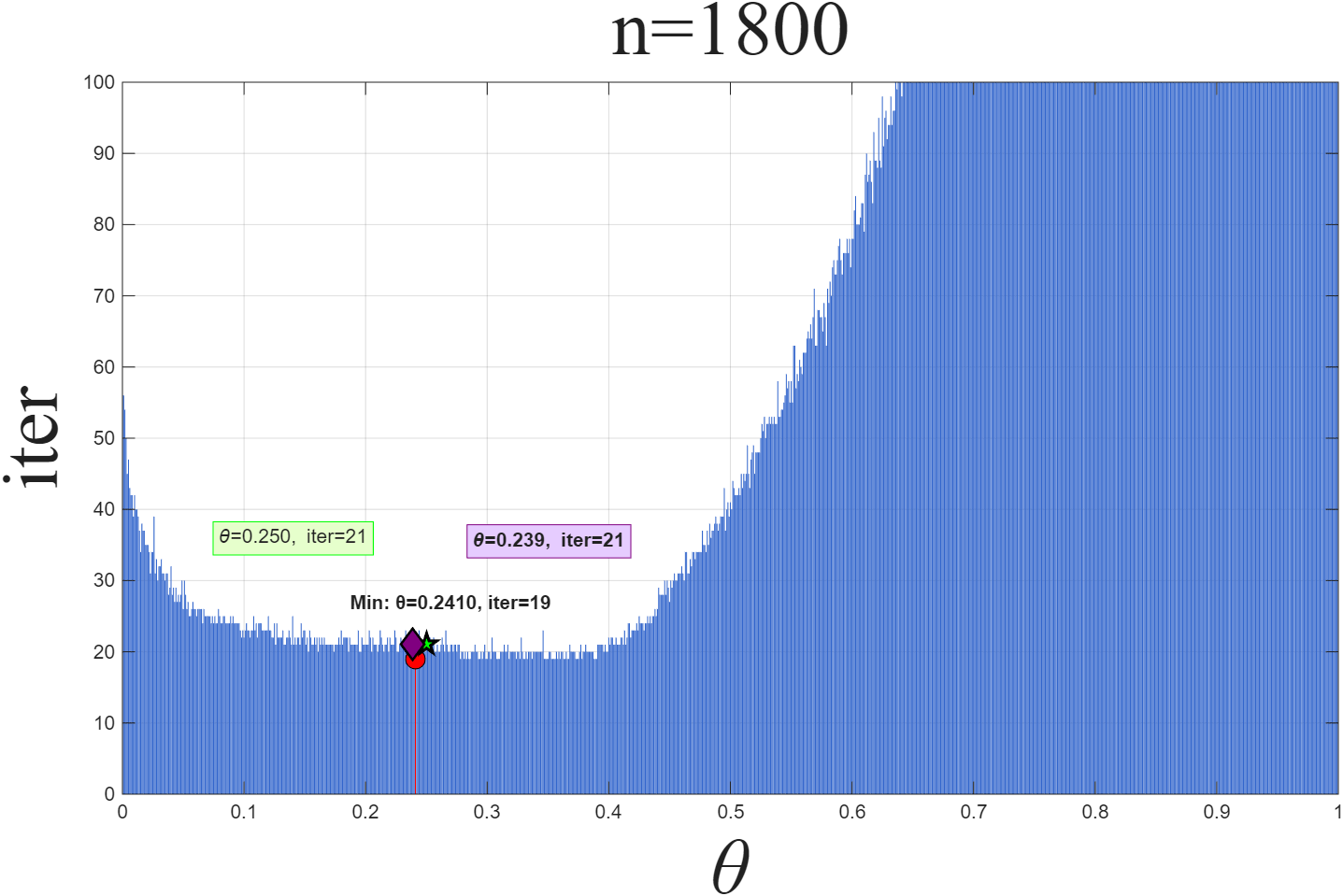}
	\end{minipage}
				\begin{minipage}{0.32\linewidth}
		\centering
		\includegraphics[width=1.1\linewidth]{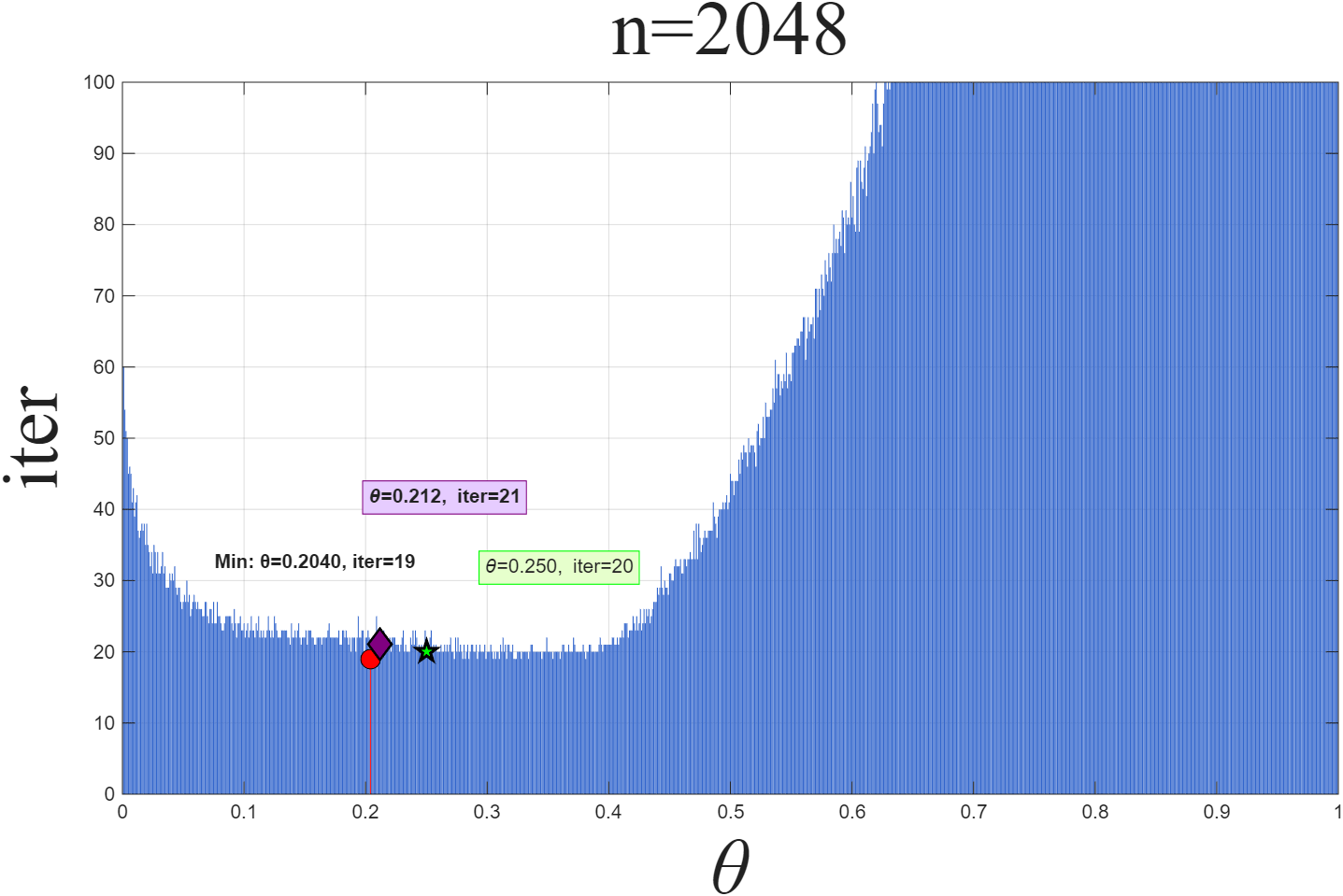}
	\end{minipage}
\end{figure}

Figure \ref{Fig:possion_pre} presents the prediction curve of { $\theta$} versus $n$ and Table \ref{tab:possion2} lists the results obtained with the traversal-optimized parameters, default parameters, and predicted parameters for {  $n=1024,~1200,~ 1400,~1600,~1800,~2048$. Figure \ref{hm_travel_more} shows the variation trend of the number of iterations with the parameter theta.} It can be observed that the number of iterations with the predicted parameters is close to that with the optimal parameters, demonstrating the effectiveness of using GPR for parameter prediction.


\begin{table}[h]
	\centering

	\caption{ For different numbers of subdivisions, the GPR predictions, the actual optimal values, the default values, and the corresponding iteration counts }
	\label{tab:possion2}
	\begin{tabular}{|c|c|c|c|c|c|c|c|c|}
		\hline
		$n$ & \makecell{predicted \\ $\theta$} & iter  & \makecell{CPU \\time}      & \makecell{optimal \\ $\theta$ }&iter &\makecell{CPU \\time} & \makecell{iter of \\ default\\ value} &   \makecell{CPU time of\\ default value}  \\ \hline
		1024  &0.368 &17  &{ 3.1793}  &0.276  & 17 & { 3.3035} &19 &{ 3.6045} \\ \hline
		1200  &0.375 &19  &{ 4.7222}   &0.22  & 18 & { 4.7302} &20 &{ 5.0649} \\ \hline
		1400  &0.35 &20 &{  6.4806}  & 0.235 &19 &{ 6.8265}   &21 &{ 7.2822} \\ \hline
		1600  &0.278  &21  &{ 9.4118}  & 0.29 &18 &{ 8.3511} & 21   &{ 9.3738}  \\ \hline
	1800 &0.239   &21 &{ 11.9171}  & 0.241 &19 &{ 11.1168} & 21  &{ 11.8335}  \\ \hline
		2048 &0.212  &21 &{ 15.9054} & 0.204 &19 &{ 14.8292} & 20 &{  15.2840} \\ \hline
	\end{tabular}
	
\end{table}

Table \ref{tab_possion_kernel1} and Table \ref{tab_possion_kernel2} provide a comparative analysis of predictive accuracy for various kernel function combinations in the context of parameter prediction for the constant coefficient Poisson equation. 
Among them, the combined Gaussian and Laplacian kernel function achieves the highest $R^2$ and Corr values.
The results indicate that the combined Gaussian and Laplacian kernel function outperforms the single Gaussian kernel across multiple evaluation metrics related to parameter prediction. Remarkably, except for the combination of rational quadratic, Gaussian kernel and  the combination of rational quadratic and Laplace kernel, which have lower 
PICP, the PICP values for other combined kernel functions are at least as high as those of the single Gaussian kernel.Among these, four combined kernel functions achieved a PICP value of 100 $\%$. This suggests that using multiple kernel functions can lead to a more robust predictive model.
Combined kernel functions significantly improve GPR model predictions while maintaining prediction reliability and potentially enhancing other performance metrics.

Combined kernel functions demonstrate superior performance in constant coefficient Poisson equation parameter prediction, providing more accurate and reliable results. This finding is crucial for developing effective GPR models, as kernel selection directly impacts predictive performance in real-world applications.

\begin{table}[H]
	\centering
	\caption{ Evaluation of kernel function combinations for the constant coefficient Poisson equation}
		\label{tab_possion_kernel1}
	\scalebox{0.65}{
		\begin{tabular}
			{|c|c|c|c|c|c|c|c|c|}
			\hline \text { Kernel function } & \text { MSE } & \text { RMSE } & \text { MAE } & R$^{ 2 } $ & \text { BIC } & \text { Corr } & \text { MdAPE } & \text { LOO-SPE } \\ \hline
			\text { Gaussian+Laplacian } &4.4889 $\times 10^{-4}$ & 0.0211 & 0.0193 & 0.3965 & -250.2432 & 0.7294 & 0.0708 & 4.6958 $\times 10^{-4}$\\ \hline
			\text { Gaussian+Exponential } & 8.0597 $\times 10^{-4}$ & 0.0283 & 0.0265 & 0.0834 & -262.6505 & 0.4165 & 0.1092 & 5.1552 $\times 10^{-4}$ \\ \hline
			\text { Gaussian } & 5.4078 $\times 10^{-4} $ & 0.0232 & 0.0179 & 0.2730 & -344.3922 & 0.6799 & 0.0421 & 4.7050 $\times 10^{-5}$
			\\ \hline
			\text { Rational Quadratic+Laplacian } & 1.8333 $\times 10^{-3}$ & 0.0428 & 0.0334 & 0.0464 & -550.7246 & 0.7820 & 0.0976 & 1.6778 $\times 10^{-4}$ \\ \hline
			\text { Mat\"{e}rn+Laplacian } & 4.6227 $\times 10^{-4}$ & 0.0215 & 0.0198 & 0.3785 & -245.4437 & 0.7207 & 0.0743 & 4.8102 $\times 10^{-4} $
			\\ \hline
			\text { Rational Quadratic+Gaussian } & 7.7147 $\times 10^{-4} $& 0.0277 & 0.0261 & 0.0370 & -259.8924 & 0.6199 & 0.1024 & 6.3889 $\times 10^{-4}$ \\ \hline
			\text { Mat\"{e}rn+Gaussian+Laplacian } & 4.6351 $\times 10^{-4}$ & 0.0215 & 0.0198 & 0.3769 & -232.9368 & 0.7208 & 0.0743 & 4.8336 $\times 10^{-4} $\\ \hline
		\end{tabular}}
\end{table}

\begin{table}[H]
	\centering
{
	\caption{PICP evaluation for the constant coefficient Poisson equation}
		\label{tab_possion_kernel2}
	\begin{tabular}{|c|c|}
		\hline
		Kernel function & PICP \\ \hline
		Gaussian+Laplacian & 100.0 \% \\ \hline
		Gaussian+Exponential & 100.0 \% \\ \hline
		Gaussian & 57.1 \% \\ \hline
		Rational Quadratic+Laplacian & 42.9 \% \\ \hline
		Mat\"{e}rn+Laplacian & 100.0 \% \\ \hline
		Rational Quadratic+Gaussian & 57.1 \% \\ \hline
		Mat\"{e}rn+Gaussian+Laplacian & 100.0 \% \\ \hline
	\end{tabular}}
\end{table}

{
\subsection{Diffusion Equation}
We consider the following two-dimensional diffusion equation, derived from \cite{Z-X-2024}:

\begin{align}
	-\nabla \cdot (\kappa \nabla u) &= f_{1}, \quad x \in \Omega, \notag \\
	u &= f_{2}, \quad x \in \partial\Omega. \label{2ks}
\end{align}
The domain is $[0,1]^{2}$, and the diffusion coefficients are
\begin{equation}
	\kappa =
	\left[
	\begin{array}{cc}
		10^{M r_{0}} & 0 \notag\\
		0 & 10^{M r_{1}} \notag 
	\end{array}
	\right],
\end{equation}
where $r_{0},r_{1}$ are random numbers in the domain $(0,1)$. The positive number $M$ will influence the multiscale properties of the matrix.

During the discretization process, the computational domain is uniformly divided into $T^{2}$ blocks based on the partition number $T$, as illustrated in Figure \ref{Fig:Block1}. The diffusion coefficient $\kappa$ remains constant within each block but  varies across different blocks. Consequently, even with fixed discretization parameter $n$ and partition number $T$, different random seeds will yield different matrices.
\begin{figure}[H]
	\centering
    \caption{When $T=2$, the computational domain is uniformly     partitioned into four blocks    $(B_{i}, i=1,2,3,4)$.     Diffusion coefficient $\kappa$ is the same in each block, when $B_{i}\neq B_{j}$, $k_{i}\neq k_{j}$. }
    \label{Fig:Block1}
		\includegraphics[width=1\linewidth]{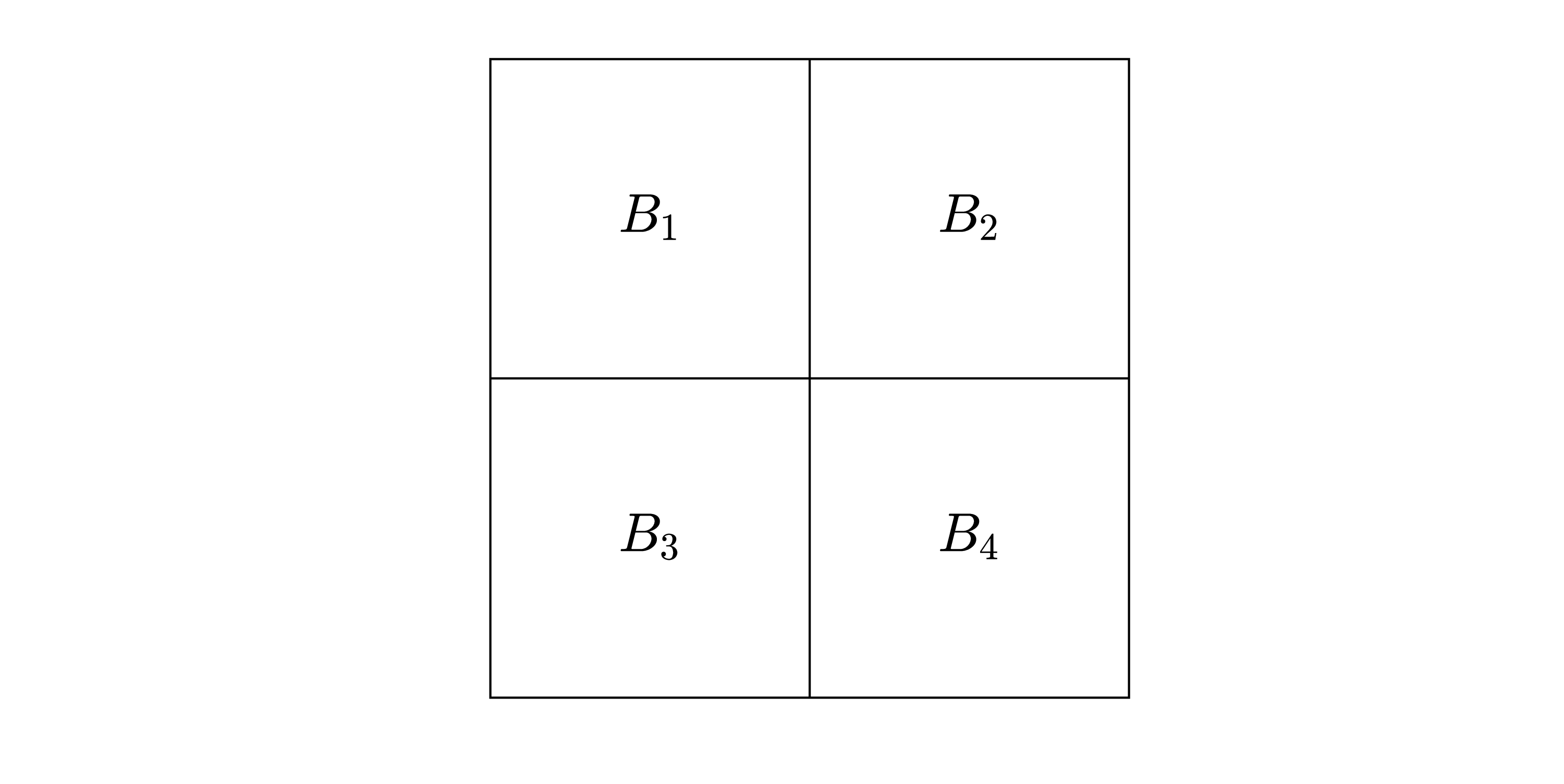}
\end{figure}

\begin{figure}[H]
	\centering
	\caption{Traversal plots for three sets of test matrices, along with iteration steps corresponding to GPR prediction, default values, and optimal parameters}
	\label{fig:diffusion_iter}
	\begin{minipage}{0.32\linewidth}
		\centering
		\includegraphics[width=1\linewidth]{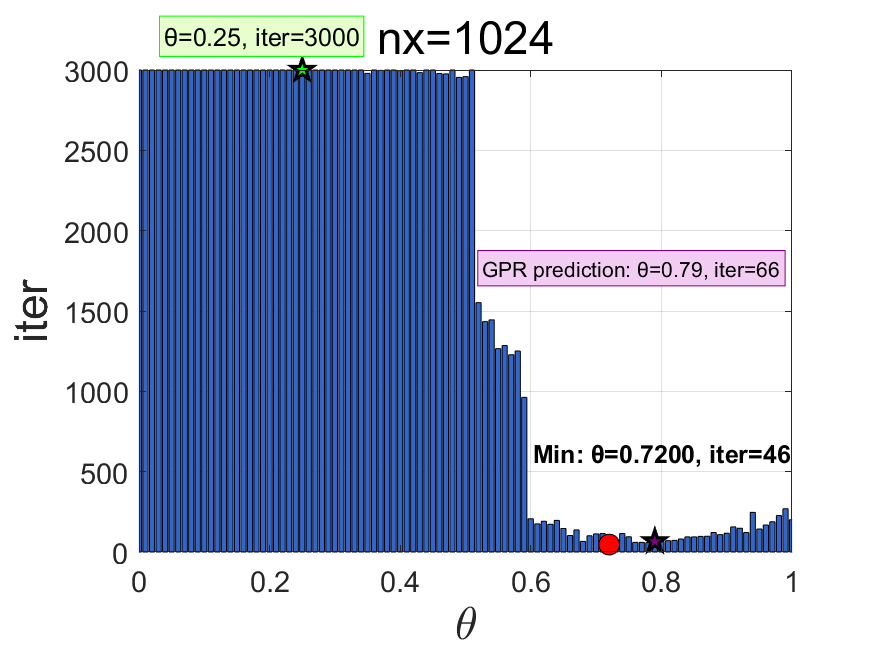}
	\end{minipage}
	\begin{minipage}{0.32\linewidth}
		\centering
		\includegraphics[width=1\linewidth]{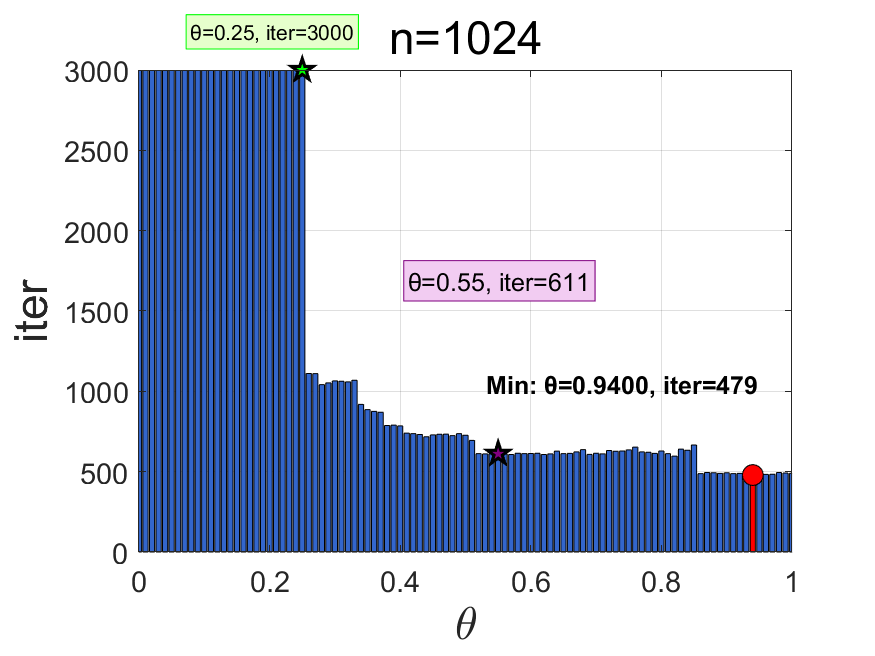}
	\end{minipage}
		\begin{minipage}{0.32\linewidth}
		\centering
		\includegraphics[width=1\linewidth]{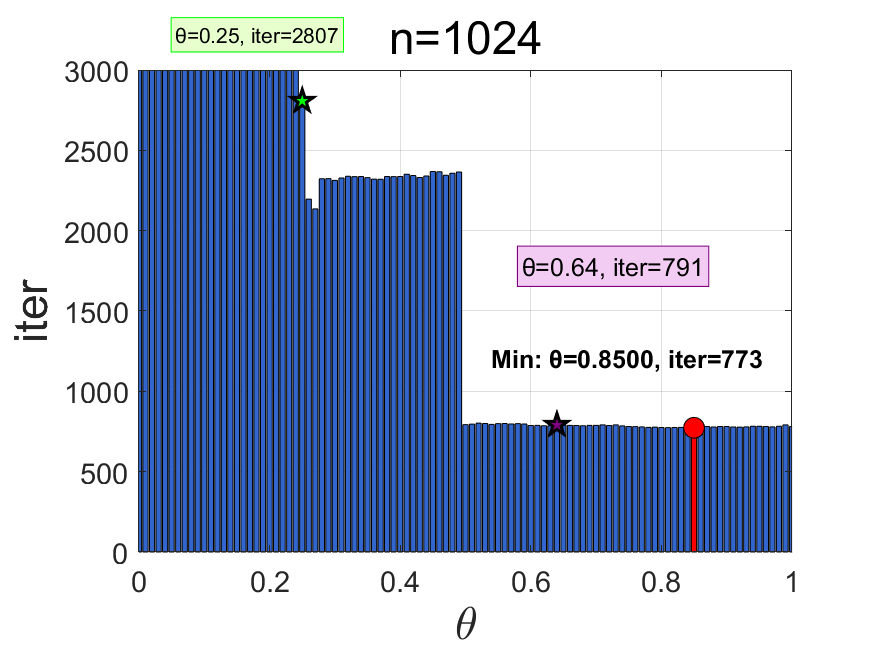}
	\end{minipage}
\end{figure}

Specifically, we generate the training set starting from $n=64$ with a step size of $16$ until $n=400$. For  each matrix generation, $T$ randomly varies within the interval $(10,20)$, and the random seeds  $ {seed}$  is equal to the index of the matrix. We will now conduct three sets of experiments. For each set, the training set consists of 23 matrices, and we will test them using matrices of size $n = 1024$. 
Figure \ref{fig:diffusion_iter} compares the iteration processes under three parameter select strategies: optimal parameters, default parameter, and GPR-predicted parameter. As shown in the Figure \ref{fig:diffusion_iter}, there is a notable difference between the default parameters and the optimal parameters in this case. Under the default parameters, the algorithm fails to converge even after 3000 iterations, whereas with the optimal parameters, convergence is achieved in fewer than 800 iterations.
To further validate the precision of GPR prediction, we test experiments using the predicted, optimal, and default $\theta$ respectively, and list the average number of iterations and CPU time from three parameter selection strategy in Table \ref{tab:diffusion2}. Results show that the GPR-predicted parameters reduce computational costs compared to default values while closely approximating the optimal parameters’ efficiency.
}

{ 
Table \ref{diffusion_PICP1} presents the evaluation of predictive performance for different kernel function combinations in the diffusion equation. According to various metrics, the Rational Quadratic+Laplacian combination achieves the best performance in MSE, MAE, and Corr. In contrast, the single Gaussian kernel yields the largest MSE. Notably, the Gaussian+Laplacian and Matérn+Laplacian combinations exhibit similar predictive performance, while hybrid combinations that include the Gaussian kernel (such as Rational Quadratic+Gaussian and Mat\'ern+Gaussian+Laplacian) show significantly higher errors.

Table \ref{diffusion_PICP2} shows PICP for each kernel combination. Rational Quadratic+Laplacian and Matérn+Laplacian achieve the highest PICP, i.e. 83.3\%, while others stay at 66.7\%. This demonstrates that hybrid kernels with Laplacian not only improve point prediction accuracy but also enhance prediction interval reliability. Overall, hybrid kernels—especially those including Laplacian—yield more accurate and robust GPR models, providing an effective kernel selection strategy for parameter inversion in complex PDEs.
}

\begin{table}[H]
\centering
	\caption{Comparison of results for AMG solving the diffusion equation}
		\label{tab:diffusion2}
		{
\begin{tabular}{|llllllll|}
\hline
\multirow{2}{*}{$n$} & \multicolumn{2}{l}{predicted $\theta$} & \multicolumn{2}{l}{optimal $\theta$} & \multicolumn{2}{l}{$\theta=0.25$} & \multirow{2}{*}{Speed up} \\ \cline{2-7}
                      & iter                    & time(s)                   & iter                   & time(s)                  & iter                  & time(s)                  &                           \\ \hline
1024 &   489.34   & 78.89   & 432.67  &  51.17           & { 2935.67}    &  178.90                &     2.27                   \\ \hline
\end{tabular}
}
\end{table}

\begin{table}[h]
\centering
{ 
\caption{Evaluation of kernel function combinations for the diffusion equation.}\label{diffusion_PICP1}
\scalebox{0.67}{ 
  \begin{tabular}{|c|c|c|c|c|c|c|c|c|}
    \hline 
    \text{Kernel function} & \text{MSE} & \text{RMSE} & \text{MAE} & $R^{2}$ & \text{BIC} & \text{Corr} & \text{MdAPE} & \text{LOO-SPE} \\
    \hline 
    \text{Gaussian+Laplacian} & $3.8986 \times 10^{-3}$ & 0.0624 & 0.0404 & 0.0433 & $-158.0800$ & 0.6646 & 0.0431 & $2.9477 \times 10^{-3}$ \\
    \hline 
    \text{Gaussian+Exponential} & $5.4413 \times 10^{-3}$ & 0.0737 & 0.0484 & 0.0352 & $-194.1815$ & 0.5201 & 0.0550 & $6.9429 \times 10^{-4}$ \\
    \hline 
    \text{Gaussian} & $6.0016 \times 10^{-3}$ & 0.0775 & 0.0534 & 0.0427 & $-218.7034$ & 0.5126 & 0.0888 & $3.8207 \times 10^{-4}$ \\
    \hline 
    \text{Rational Quadratic+Laplacian} & $3.7333 \times 10^{-3}$ & 0.0611 & 0.0406 & 0.0838 & $-155.2654$ & 0.6993 & 0.0499 & $2.0250 \times 10^{-3}$ \\
    \hline 
    \text{Mat\'ern+Laplacian} & $3.9262 \times 10^{-3}$ & 0.0627 & 0.0405 & 0.0365 & $-145.8736$ & 0.6608 & 0.0436 & $2.1008 \times 10^{-3}$ \\
    \hline 
    \text{Rational Quadratic+Gaussian} & $6.6051 \times 10^{-3}$ & 0.0812 & 0.0528 & 0.0621 & $-205.9296$ & 0.4080 & 0.0583 & $4.6941 \times 10^{-4}$ \\
    \hline 
    \text{Mat\'ern+Gaussian+Laplacian} & $6.7545 \times 10^{-3}$ & 0.0822 & 0.0531 & 0.0657 & $-186.6553$ & 0.3913 & 0.0567 & $5.4013 \times 10^{-4}$ \\
    \hline
  \end{tabular}
}
}
\end{table}

\begin{table}[h]
\centering 
{ 
\caption{\raggedright  PICP evaluation for the diffusion equation.}\label{diffusion_PICP2}
\begin{tabular}{|c|c|}
\hline \text { Kernel function
 } & \text { PICP }  \\
\hline \text { Gaussian+Laplacian
 } & 66.7 \% \\
\hline \text { Gaussian+Exponential
 } & 66.7 \% \\
\hline \text { Gaussian } & 66.7 \% \\
\hline \text { Rational Quadratic+Laplacian
 } & 83.3 \% \\
\hline \text { Matérn+Laplacian
 } & 83.3 \% \\
\hline \text { Rational Quadratic+Gaussian
 } & 66.7 \%  \\
\hline \text { Matérn+Gaussian+Laplacian
 } & 66.7 \%  \\
\hline
\end{tabular}
}
\end{table}

 { In \cite{Z-X-2024 } the authors introduce  an AutoAMG ($\theta$) method for predicting the strong threshold parameter in AMG based on graph neural networks and present the following experiment on 2D diffusion equations: a training set of 80 matrices and a test set of 20 matrices are constructed, with partition number $n$ randomly selected from $(50, 100)$ and $T$ from $(10, 20)$.
The following Table \ref{db1} presents the experimental results obtained from \cite{Z-X-2024 } , and Table \ref{db2} displays our experimental results. Please note that here ``nrow" denotes the order of the test matrix, and the data for each metric are averaged over 20 experiments.

\begin{table}[h]
\centering
{ 
\caption{Test results of AutoAMG($\theta$).}\label{db1}
\begin{tabular}{cccccccc}
\toprule
\multirow{2}{*}{nrow} & \multicolumn{2}{c}{optimal $\theta$} & \multicolumn{2}{c}{$\theta = 0.25$} & \multicolumn{2}{c}{AutoAMG($\theta$)} & \multirow{2}{*}{speedup} \\
\cmidrule(lr){2-3} \cmidrule(lr){4-5} \cmidrule(lr){6-7}
& iter & time(s) & iter & time(s) & iter & time(s) & \\
\midrule
5659 & 185.25 & 0.15 & 496.20 & 0.38 & 257.30 & 0.21 & 1.81 \\
\bottomrule
\end{tabular}}
\end{table}

\begin{table}[h]
\centering
{ 
\caption{Test results of GPR.}\label{db2}
\begin{tabular}{cccccccc}
\toprule
\multirow{2}{*}{nrow} & \multicolumn{2}{c}{optimal $\theta$} & \multicolumn{2}{c}{$\theta = 0.25$} & \multicolumn{2}{c}{GPR predicted $\theta$} & \multirow{2}{*}{speedup} \\
\cmidrule(lr){2-3} \cmidrule(lr){4-5} \cmidrule(lr){6-7}
& iter & time(s) & iter & time(s) & iter & time(s) & \\
\midrule
5659 & 99.70 & 0.12 & 709.75 & 0.61 & 131.50 & 0.19 & 3.21 \\
\bottomrule
\end{tabular}}
\end{table}

It should be noted that the samples tested in our experiment and those in \cite{Z-X-2024 }  may not be identical (as the test matrices were randomly generated). However, the experimental data, to some extent, demonstrates the comparable effectiveness of our GPR to Autoamg($\theta$) and highlights its significant research implications in the field of data-driven AMG parameter tuning.}

\subsection{{Parabolic Equation}}

Consider the following parabolic equation:

$$\begin{cases}
\displaystyle
\frac{\partial u}{\partial t} - \nabla \!\cdot\! \bigl(\kappa \nabla u\bigr) = f_{1},
& (x,y)\in\Omega,\; t\in(0,1],\\[4pt]
u(x,y,0) = \cos(\pi x)*\cos(\pi y), & (x,y)\in\Omega,\\[4pt]
u(x,y,t) = f_{3}, & (x,y)\in\partial\Omega,\; t\in[0,1].
\end{cases}$$

The exact solution to this equation is $u=\cos(\pi x)\cos(\pi y) e^{-\frac{\pi^{2}}{8}t} $,  where $\kappa$ is obtained from (\ref{2ks}), and $f_{3}$   is derived from $f_{2}$   in (\ref{2ks}). Discretization is performed using the finite element method, with a time step size of $\Delta t=0.25$.  For a fixed mesh partition number n, discretization in the time direction yields four discrete time instances, resulting in four linear systems. Our focus remains on the linear system at the final time instant $t_{n}=1$. The traversal for the $\theta$  covers the interval $[0,1]$ with a step size of $0.01$, and the iteration process is limited to a maximum of 1500 iterations. The training set is constituted by selecting   $n$ from 64 to 400, incrementing by $\Delta n=16$.  The detailed situation is listed in Table \ref{tab:para_traversal}.

\begin{table}[h]
	\centering
	\caption{Traversal results of AMG solving the parabolic equation }
	\label{tab:para_traversal}
	\begin{tabular}{|c|c|c|c|c|c|}
		\hline
		$n$ & $\theta$ & iter & $n$ & $\theta$ & iter \\ \hline
		64  &0.83 &197 & 240& 0.35  &296 \\ \hline
		80  &0.41 &70 & 256 &0.9  &197 \\ \hline
		96  &0.99  &79 & 272 &0.4 & 226  \\ \hline
		112 &0.64  &28 & 288 &0.34 & 410 \\ \hline
		128 &0.98  &110 & 304 &0.36 & 308 \\ \hline
		144 &0.45  &239  & 320 &0.63 & 27 \\ \hline
		160 &0.98  &97 & 336& 0.64 & 243 \\ \hline
		176 &0.99  &132 &352& 0.57 & 328 \\ \hline
	    192 &0.98  &165& 368 &0.39 & 402 \\ \hline
		208 &0.62  &218 & 384 &0.89 &243 \\ \hline
		224 &0.88  &330 & 400 &0.95 & { 255} \\ \hline
	\end{tabular}
\end{table}

\begin{figure}[H]
	\centering
	\caption{Traversal plots for three sets of test matrices, along with iteration steps corresponding to GPR prediction, default values, and optimal parameters}
	\label{fig_para_pre}
	\begin{minipage}{0.32\linewidth}
		\centering
		\includegraphics[width=1\linewidth]{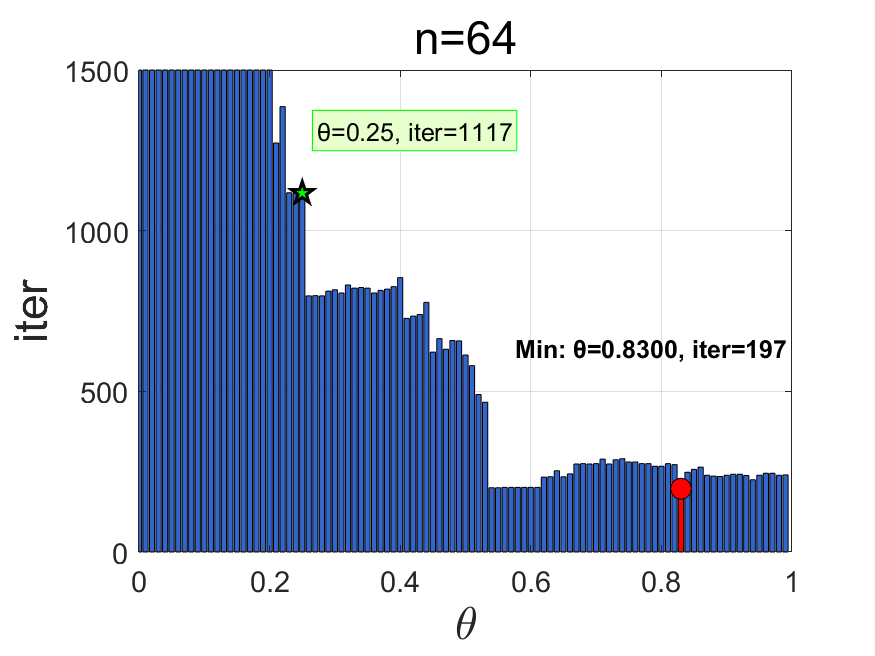}
	\end{minipage}
	\begin{minipage}{0.32\linewidth}
		\centering
		\includegraphics[width=1\linewidth]{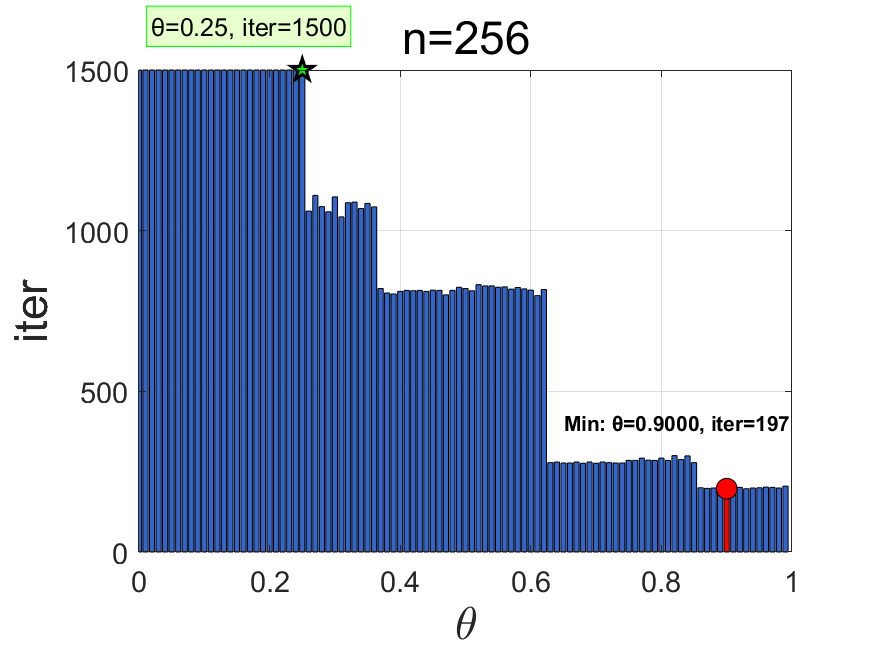}
	\end{minipage}
		\begin{minipage}{0.32\linewidth}
		\centering
		\includegraphics[width=1\linewidth]{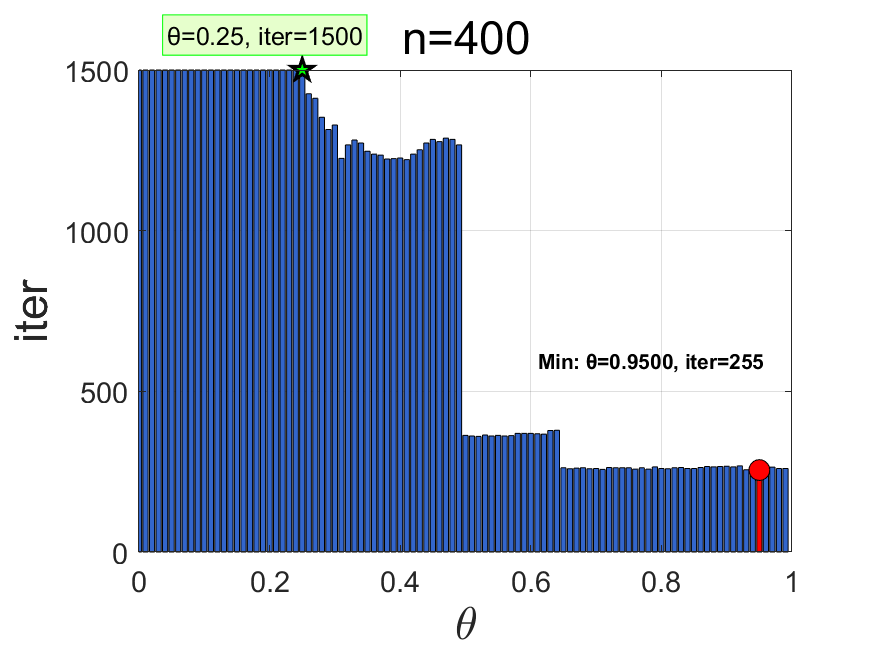}
	\end{minipage}
\end{figure}

Figure \ref{fig_para_pre} illustrates the traversal scenarios for $n=64,~256,~400$. As can be seen from the figure, when $\theta = 0.25$ (the default parameter), AMG fails to converge after 1500 iterations,and there exists a significant difference between the optimal and default parameters.

\begin{figure}[H] 
	\centering
	\caption{Regression curves of $\theta$ with respect to $n$ predicted by GPR for AMG solving the parabolic equation }
	\label{parabolic GPR}
	\begin{minipage}{0.32\linewidth}
		\centering
		\includegraphics[width=0.9\linewidth]{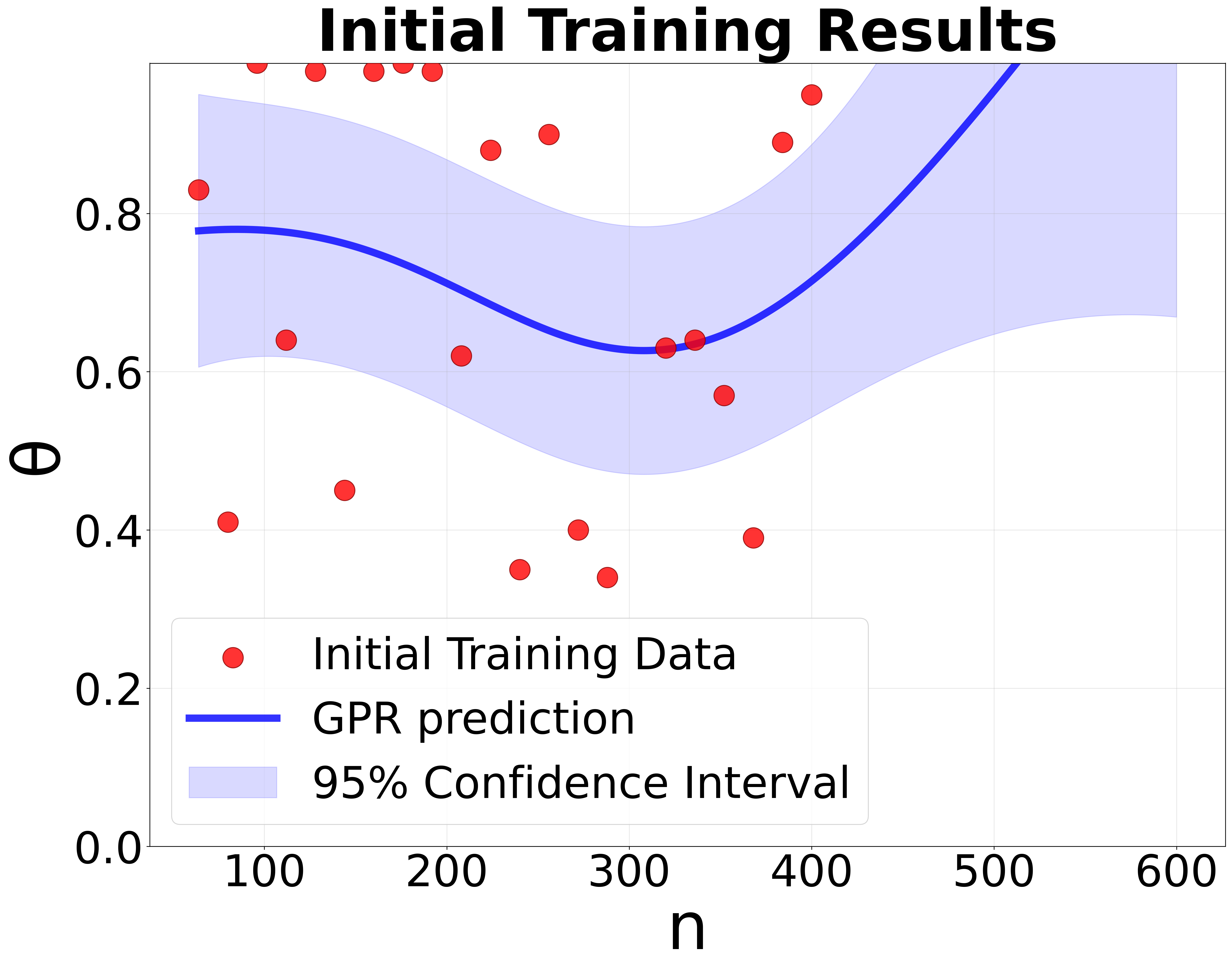}
	\end{minipage}
	\begin{minipage}{0.32\linewidth}
		\centering
		\includegraphics[width=0.9\linewidth]{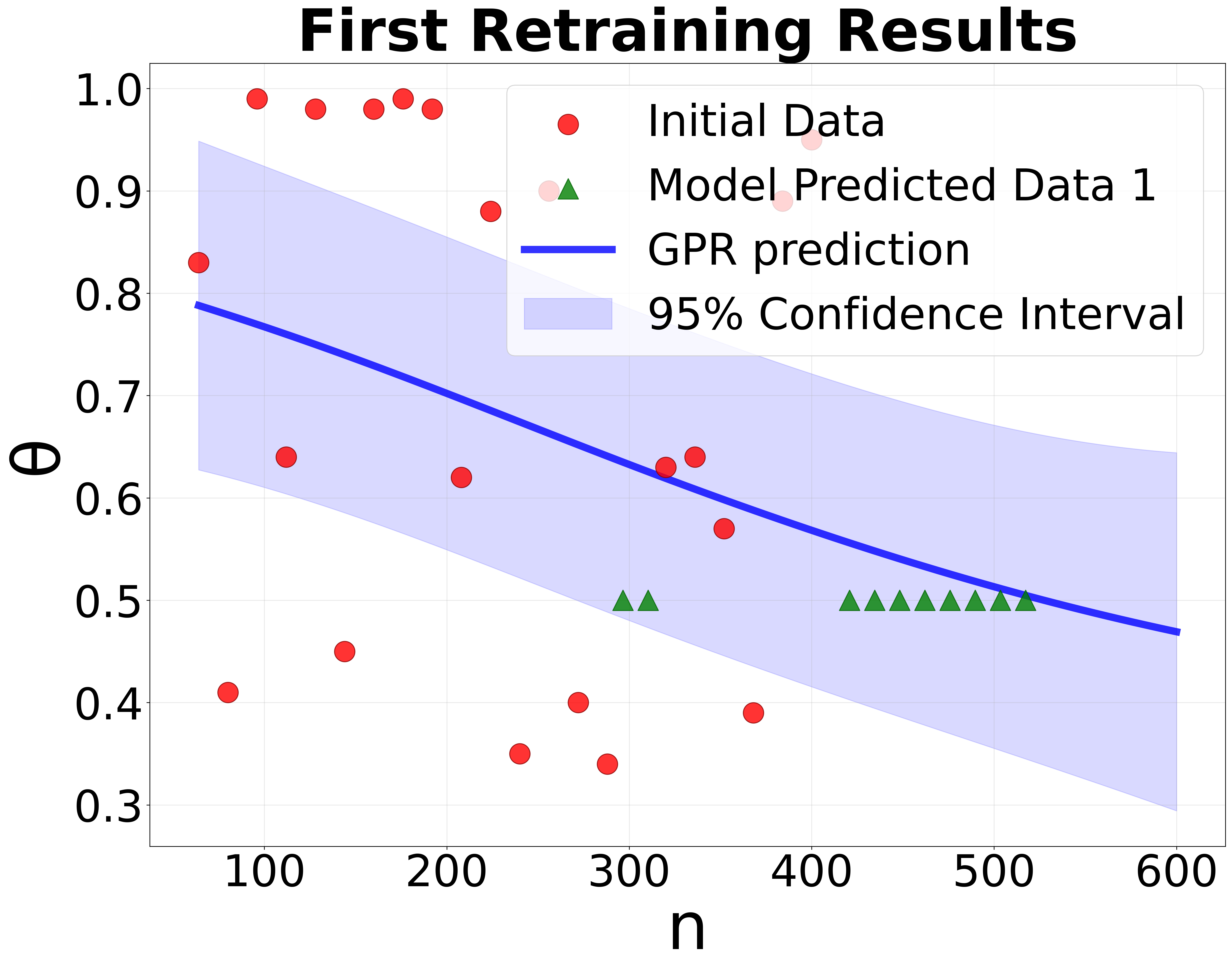}
	\end{minipage}
		\begin{minipage}{0.32\linewidth}
		\centering
		\includegraphics[width=0.9\linewidth]{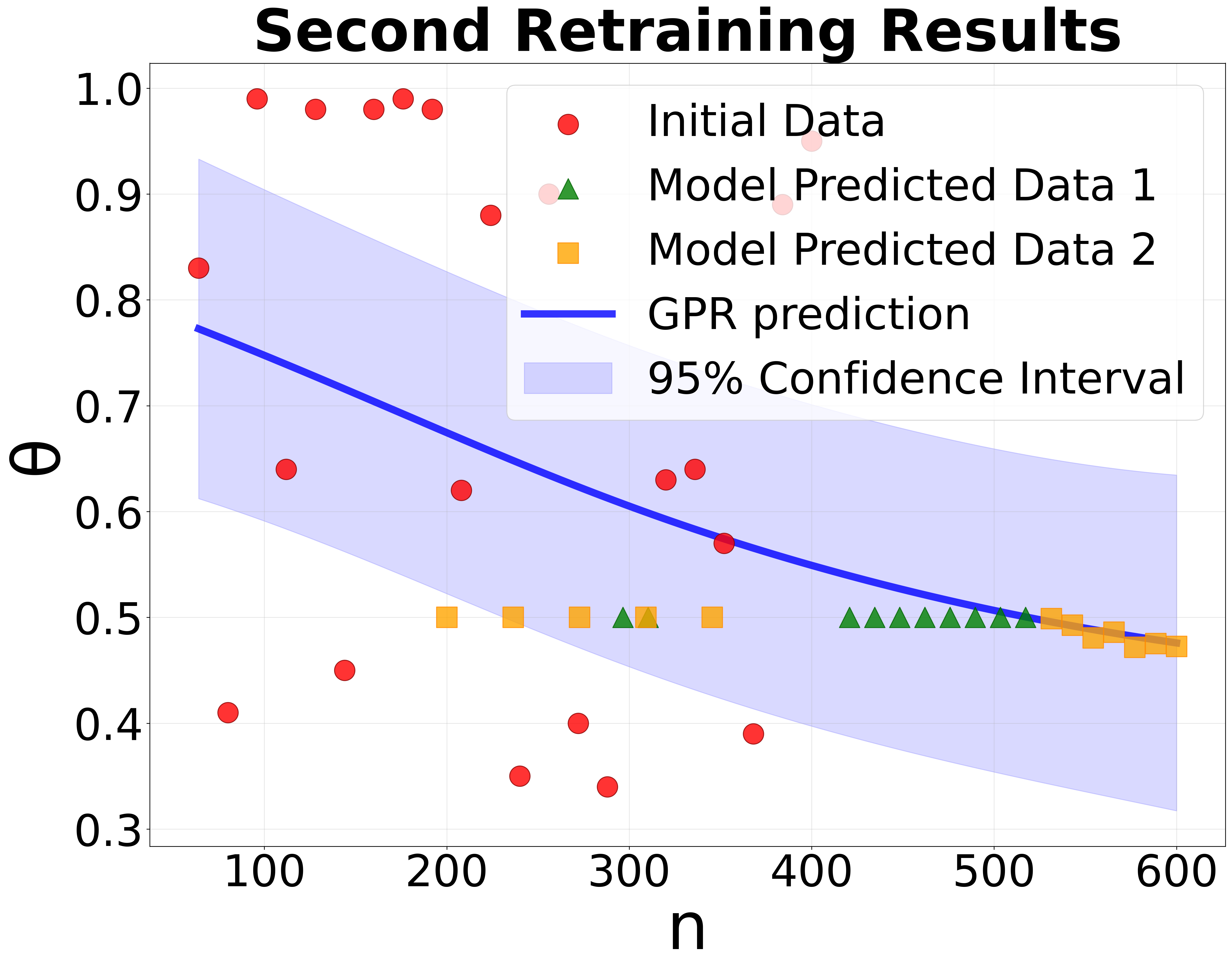}
	\end{minipage}
\end{figure}

{For retraining we adopt targeted sampling: the first retraining phase selects 10 new points uniformly distributed in [200,600], followed by a second retraining phase where 12 additional points are chosen from [200,600] to refine the model further. We summarize this in Table \ref{para1}. Figure \ref{parabolic GPR} illustrates the fitting curve of the GPR predicted values.       } 

\begin{table}[H] 
	\centering
{
\caption{Training set, and retraining set for using GPR to predict the optimal parameters for AMG solving the parabolic equation}
\label{para1}
	\scalebox{1}{
		\begin{tabular}{|c|c|c|}
			\hline
			Training set & $n\in[64, 400],~\Delta n=16$\\
			Retraining set 1 & $n\in[200,600]$, 10 randomly selected points\\
			Retraining set 2 & $n\in[200,600]$,12 randomly selected points \\
			\hline
		\end{tabular}
	}
}
\end{table}

Table \ref{para_tab_pre} lists the number of iterations and CPU time of AMG under three different parameters for $n = 1024$. AMG converges successfully with both predicted parameter $\theta$ obtained by GPR and the optimal parameter, exhibiting the number of iterations below 200 in both cases, lower than the number of iterations of AMG under the default parameter. Table \ref{para_tab_pre} shows that at 1500 iteration steps does not converge. Besides, at this point, the relative residual is $10^{-7}$, which fails to meet the preset threshold of $10^{-8}$. This demonstrates the effectiveness of the provided predictive parameters.

\begin{table}[H]
\centering
{
\caption{Comparison of results for AMG solving the Parabolic Equation}
\label{para_tab_pre}
\begin{tabular}{|lllllll|}
\hline
\multirow{2}{*}{$n$} & \multicolumn{2}{l}{predicted $\theta=0.96$} & \multicolumn{2}{l}{optimal $\theta=0.87$} & \multicolumn{2}{l|}{$\theta=0.25$}  \\ \cline{2-7}
    & iter   & time(s)     & iter      & time(s) & iter   & time(s)                               \\ \hline
1024 &   189  &100.12           &  120       & 76.12           & $>1500$    & -                            \\ \hline
\end{tabular}
}
\end{table}

{ 
Table \ref{para:PICP1} presents the predictive performance of different kernel combinations for the parabolic equation. Among the evaluated metrics, the Rational Quadratic+Gaussian combination delivers the lowest MSE, RMSE, and MAE, along with the highest $R^2$ and Corr, indicating the best point prediction accuracy.  Combinations incorporating Laplacian (such as Gaussian+Laplacian and Matérn+Gaussian+Laplacian) generally achieve lower BIC values, suggesting better model fit after accounting for complexity.
Table \ref{para:PICP2} reports the PICP values for each kernel combination. Most combinations reach a PICP of 83.3\%, demonstrating reliable prediction intervals, except for Rational Quadratic+Laplacian, which is lower at 66.7\%. Overall, hybrid kernels—particularly those combining Gaussian with other bases—provide superior accuracy and robust uncertainty quantification, offering a practical strategy for kernel selection in parameter inversion tasks involving parabolic PDEs.
}

\begin{table*}[h]
\centering
{ 
\caption{ Evaluation of kernel function combinations for the  parabolic equation.}\label{para:PICP1}%
\scalebox{0.67}{ 
  \begin{tabular}
    {|c|c|c|c|c|c|c|c|c|}

\hline \text { Kernel function } & \text { MSE } & \text { RMSE } & \text { MAE } & $R^{ 2 }$  & \text { BIC } & \text { Corr } & \text { MdAPE } & \text { LOO-SPE } \\
\hline \text { Gaussian+Laplacian } & $8.1900 \times 10^{-4}$ & 0.0286 & 0.0154 & 0.8855 & -162.0790 & 0.9646 & 0.0632 & $1.4711 \times 10^{-2}$ \\
\hline \text { Gaussian+Exponential } & $8.4267 \times 10^{-4}$ & 0.0290 & 0.0154 & 0.8851 & -182.3871 & 0.9646 & 0.0434 & $1.3765 \times 10^{-2}$ \\
\hline \text { Gaussian } & $2.1115 \times 10^{-3}$ & 0.0459 & 0.0222 & 0.8626 & -118.8061 & 0.9582 & 0.0487 & $4.2660 \times 10^{-3}$ \\
\hline \text { Rational Quadratic+Laplacian } & $8.4113 \times 10^{-4}$ & 0.0290 & 0.0148 & 0.8851 &-246.7741 & 0.9651 & 0.0846 & $8.8024 \times 10^{-1}$ \\
\hline \text { Mat\'ern+Laplacian } & $8.7103 \times 10^{-3}$ & 0.0933 & 0.0456 & 0.8460 & -158.3641 & 0.9076 & 0.0534 & $1.2143 \times 10^{-2}$ \\
\hline \text { Rational Quadratic+Gaussian } & $6.1593 \times 10^{-4}$ & 0.0248 & 0.0129 & 0.8891 & -120.4146 & 0.9663 & 0.0612 & $2.3037 \times 10^{-2}$ \\
\hline \text { Mat\'ern+Gaussian+Laplacian } & $8.1128 \times 10^{-4}$ & 0.0284 & 0.0151 & 0.8856 & -156.6667 & 0.9649 & 0.0524 & $1.3468 \times 10^{-2}$ \\
\hline
  \end{tabular}
}
}
\end{table*}

\begin{table*}[h]
 
\caption{\raggedright  PICP evaluation for the  parabolic equation.}\label{para:PICP2}%
\centering 
\begin{tabular}{|c|c|}
\hline \text { Kernel function
 } & \text { PICP }  \\
\hline \text { Gaussian+Laplacian
 } & 83.3 \% \\
\hline \text { Gaussian+Exponential
 } & 83.3 \% \\
\hline \text { Gaussian } & 83.3 \% \\
\hline \text { Rational Quadratic+Laplacian
 } & 66.7 \% \\
\hline \text { Matérn+Laplacian
 } & 83.3 \% \\
\hline \text { Rational Quadratic+Gaussian
 } & 83.3 \%  \\
\hline \text { Matérn+Gaussian+Laplacian
 } & 83.3 \%  \\
\hline
\end{tabular}

\end{table*}

\subsection{Helmholtz Equation}

We have delved into the two-dimensional Helmholtz equation, a milestone partial differential equation in physics and engineering. The Helmholtz equation is widely recognized for its significant role in modeling wave propagation phenomena, as it can describe the propagation characteristics of waves in various media, including sound waves, light waves, and other types of waves. This equation is not only central to theoretical research but also extremely important in practical applications, such as in acoustics, optics, electromagnetism, and quantum mechanics. The mathematical expression of this equation is as follows:

\[
\begin{cases}
	\Delta u + k^2 u = f, & (x,y) \in \Omega = (-1,1)^2, \\
	u(-1,y) = u(1,y) = 0, & y \in (-1,1), \\
	u(x,-1) = u(x,1) = -\sin(\pi x), & x \in (-1,1).
\end{cases}
\]
The exact solution to this equation is \( u = \sin(\pi x) \cos(\pi y) \), and the corresponding right-hand side term is \( f = (k^2 - \pi^2) \sin(\pi x) \cos(\pi y) \).

In this particular scenario, we assign the wave number \( k \) a value of {\( 2\pi \)}. The connectivity parameter is systematically traversed with an increment of 0.001 across the interval {\([0,1]\)}, while the number of partitions \(n \) spans from {64 to 400}, incrementing by 16 for each step. The detailed traversal outcomes are presented in Table \ref{Tab:hm_traversal}.

\begin{table}[h]
	\centering
{
	\caption{Traversal results of AMG solving the Helmholtz equation with coefficient {$2\pi$}}
	\label{Tab:hm_traversal}
	\begin{tabular}{|c|c|c|c|c|c|}
		\hline
		$n$ & $\theta$ & iter & $n$ & $\theta$ & iter \\ \hline
		64  &0.255 &17 & 240& 0.139  &20 \\ \hline
		80  &0.414 &19 & 256 &0.121  &21 \\ \hline
		96  &0.108  &20 & 272 &0.142 & 21  \\ \hline
		112 &0.083  &22 & 288 &0.094 & 23  \\ \hline
		128 &0.032  &31 & 304 &0.096 & 23 \\ \hline
		144 &0.260  &21  & 320 &0.066 & 24 \\ \hline
		160 &0.239  &20 & 336& 0.068 & 23 \\ \hline
		176 &0.204  &20 &352& 0.275 & 22 \\ \hline
	    192 &0.250  &20& 368 &0.328 & 22 \\ \hline
		208 &0.155  &20 & 384 &0.344 &22 \\ \hline
		224 &0.160  &20 & 400 &0.266 & 23 \\ \hline
	\end{tabular}
	}
\end{table}

\begin{figure}[H]
	\centering
	\caption{Variation of the number of iterations for AMG solving the Helmholtz equation with coefficient {$2\pi$} and different connectivity parameters $\theta$ when $n$ is {64,256,400}, respectively}
	\label{fig:hm_bianli}
	\begin{minipage}{0.325\linewidth}
		\centering
		\includegraphics[width=1\linewidth]{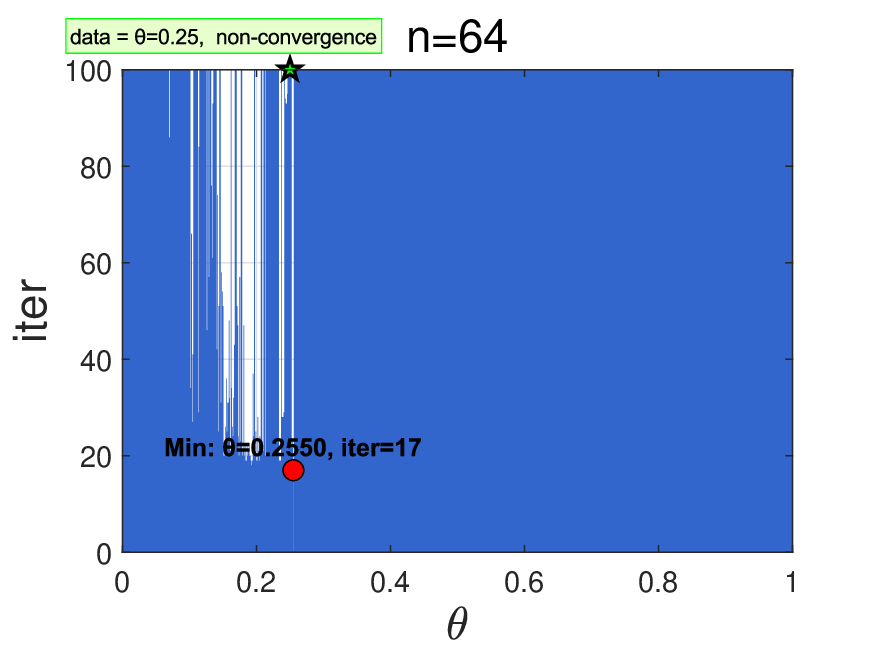}
	\end{minipage}
	\begin{minipage}{0.325\linewidth}
		\centering
		\includegraphics[width=1\linewidth]{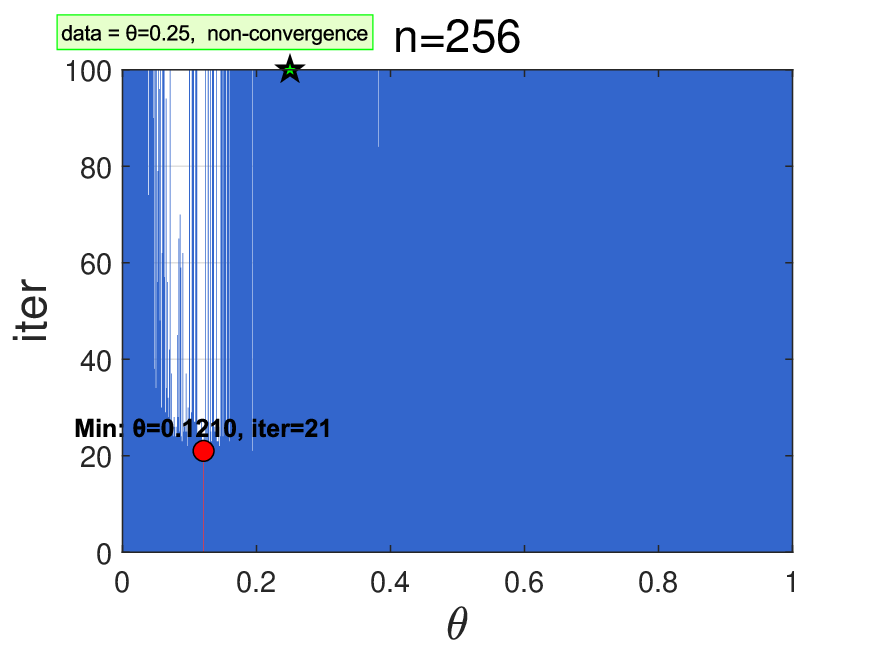}
	\end{minipage}
	\begin{minipage}{0.325\linewidth}
		\centering
		\includegraphics[width=1\linewidth]{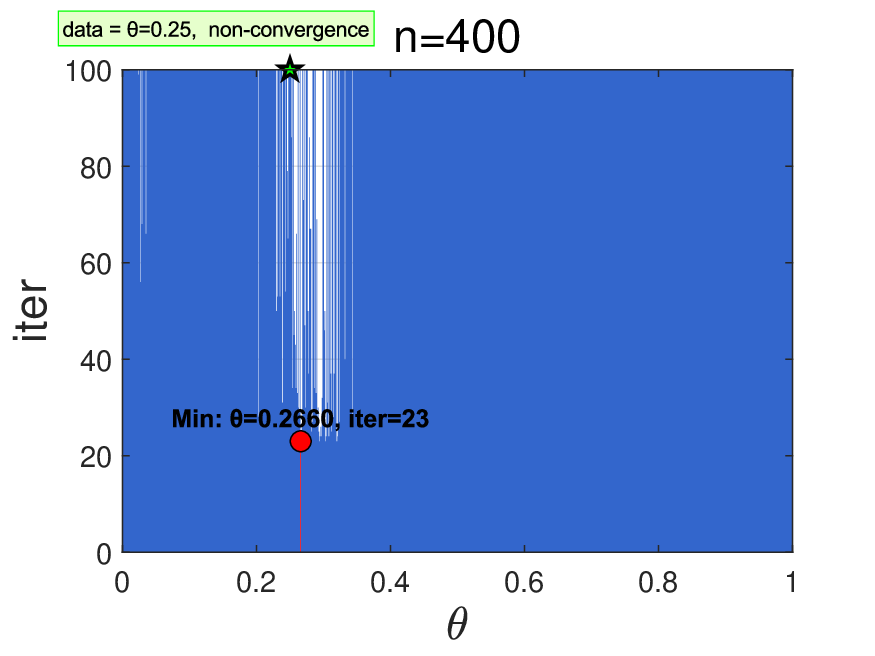}
	\end{minipage}
\end{figure}

As shown in Figure \ref{fig:hm_bianli}, we present the variation in the number of iterations for the AMG method as the traversal parameters are altered, with the number of subdivisions \(n \) specified as {64,~256,~400}. The connectivity parameter is traversed over the domain {\((0,1]\)} with a step size of 0.001, and the iteration process is limited to a maximum of 100 iterations. { In this 2d helmholtz equation, with the default parameter $\theta = 0.25$, AMG does not converge under default parameters.}

Below, we forecast the optimal parameters for the AMG method when solving the Helmholtz equation using GPR method. Table \ref{tab:hm_xunlianji} presented details of the training set, and retraining set for GPR. The training set is derived by selecting \(n \) values from {64 to 400}, incrementing by \(\Delta n = 16\).
{For retraining test, we adopt targeted sampling: the first retraining phase selects 10 new points uniformly distributed in [200,600], followed by a second retraining phase where 12 additional points are chosen from [200,600] to refine the model further. We   summarize it in Table \ref{tab:hm_xunlianji}. }
\begin{table}[H]
	\centering
\caption{Training set, and retraining set for using GPR to predict the optimal parameters for AMG solving the Helmholtz equation with coefficient {$2\pi$}}
\label{tab:hm_xunlianji}
	\scalebox{1}{
		\begin{tabular}{|c|c|c|}
			\hline
			Training set & $n\in[64, 400],~\Delta n=16$\\
			Retraining set 1 & $n\in[200,600]$, 10 randomly selected points\\
			Retraining set 2 & $n\in[200,600]$,12 randomly selected points \\
			\hline
		\end{tabular}
	}
\end{table}

\begin{figure}[H]
	\centering
	\caption{Regression curves of $\theta$ with respect to $n$ predicted by GPR for AMG solving the Helmholtz equation with coefficient $2\pi$}
	\label{fig:hm:pre}
	\begin{minipage}{0.32\linewidth}
		\centering
		\includegraphics[width=0.9\linewidth]{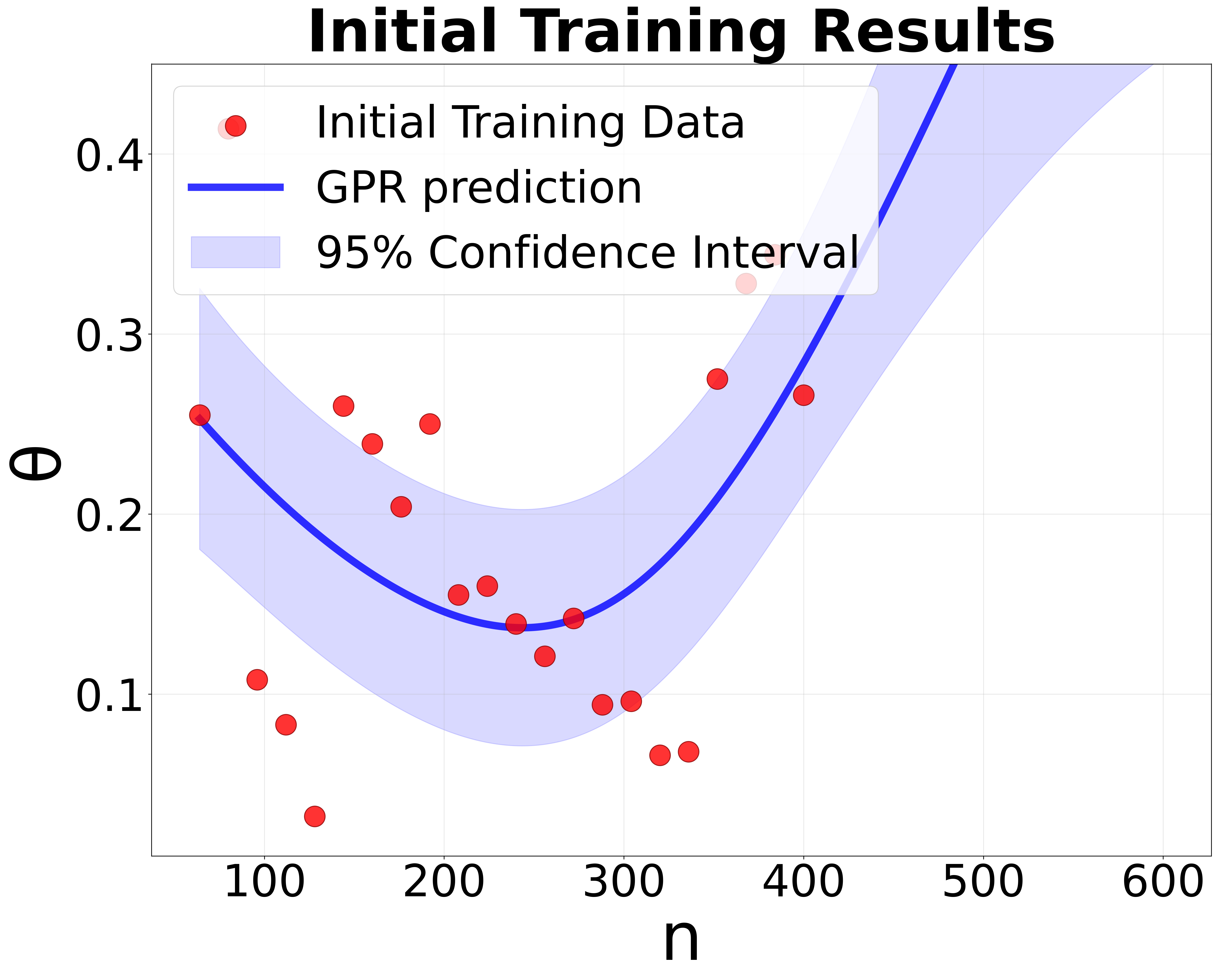}
	\end{minipage}
	\begin{minipage}{0.32\linewidth}
		\centering
		\includegraphics[width=0.9\linewidth]{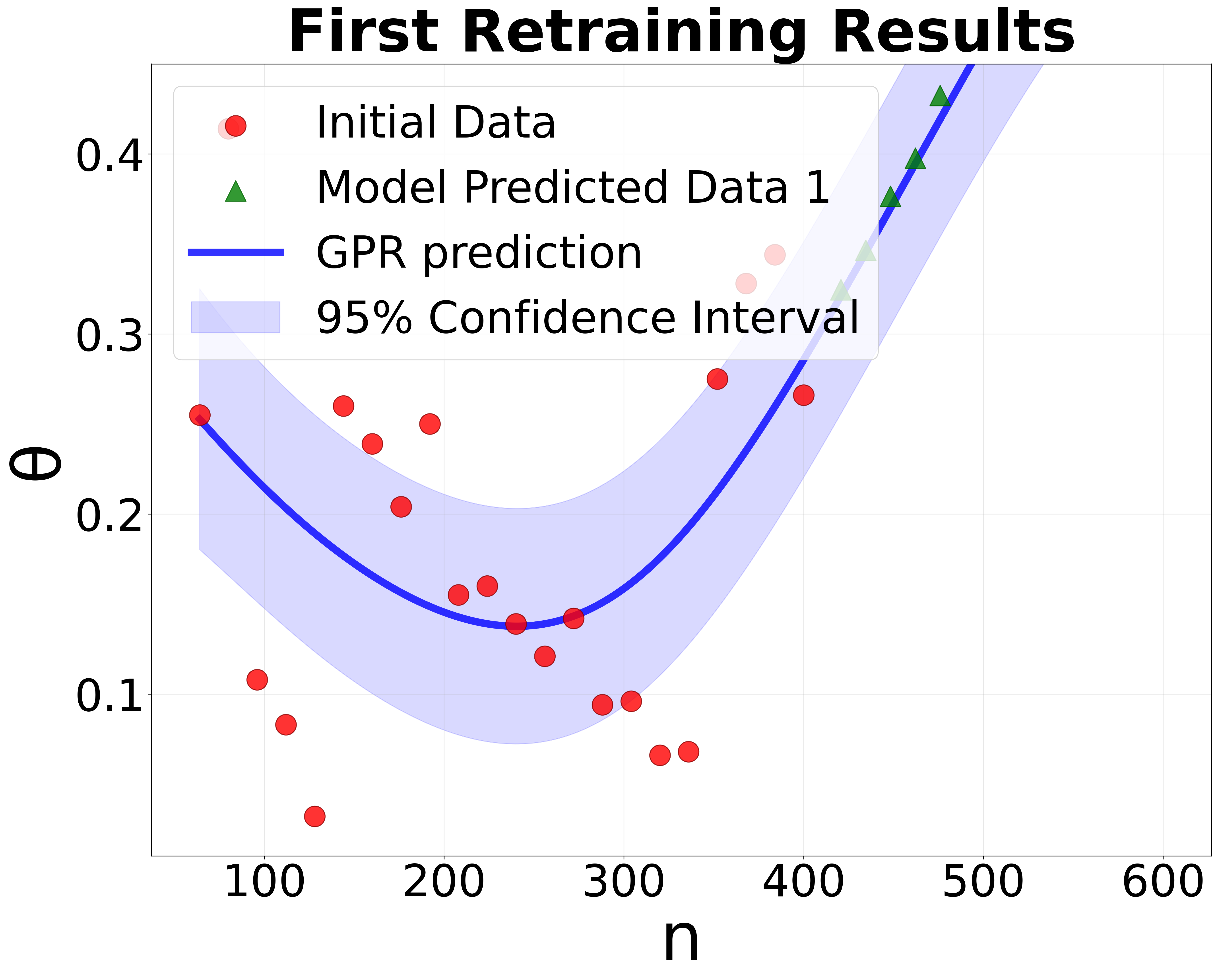}
	\end{minipage}
		\begin{minipage}{0.32\linewidth}
		\centering
		\includegraphics[width=0.9\linewidth]{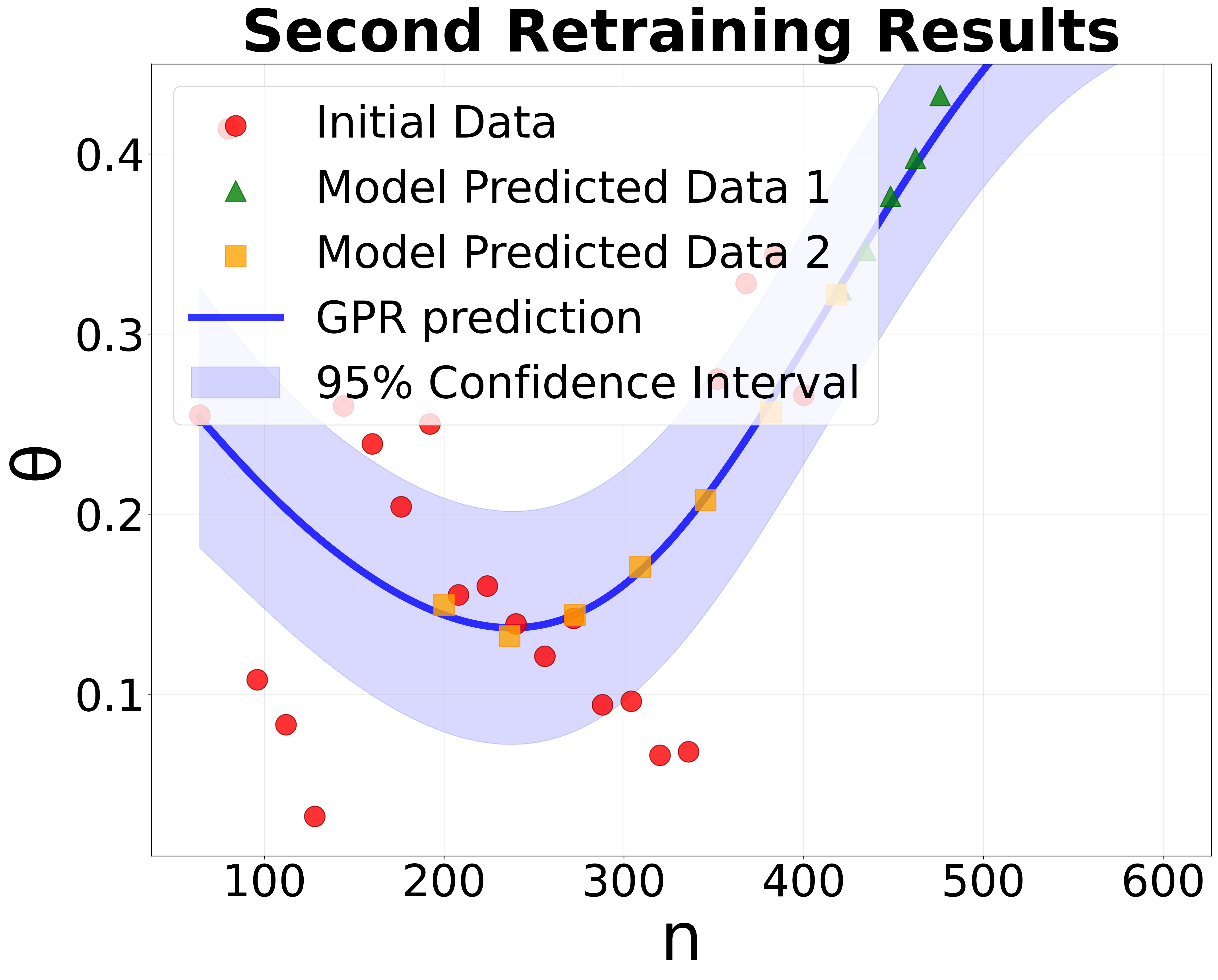}
	\end{minipage}
\end{figure}

Continuing, we depict the regression curve that illustrates the relationship between the connectivity parameter \(\theta\) and the number of partitions \(n \) in the AMG method as predicted by GPR. The blue regions within the Figure \ref{fig:hm:pre} denote the confidence intervals, signifying that the predicted results have a significant degree of reliability.

\begin{figure}[H]
	\centering
	\caption{{  Variation of the number of iterations for AMG solving the Helmholtz equation with coefficient $2\pi$ and different connectivity parameters $\theta$ when $n$ is~ 1024,~1200,~1400,~1600,~1800,~2048, respectively}}
	\label{hm_travel_more_Helm}
		\begin{minipage}{0.32\linewidth}
		\centering
		\includegraphics[width=1.1\linewidth]{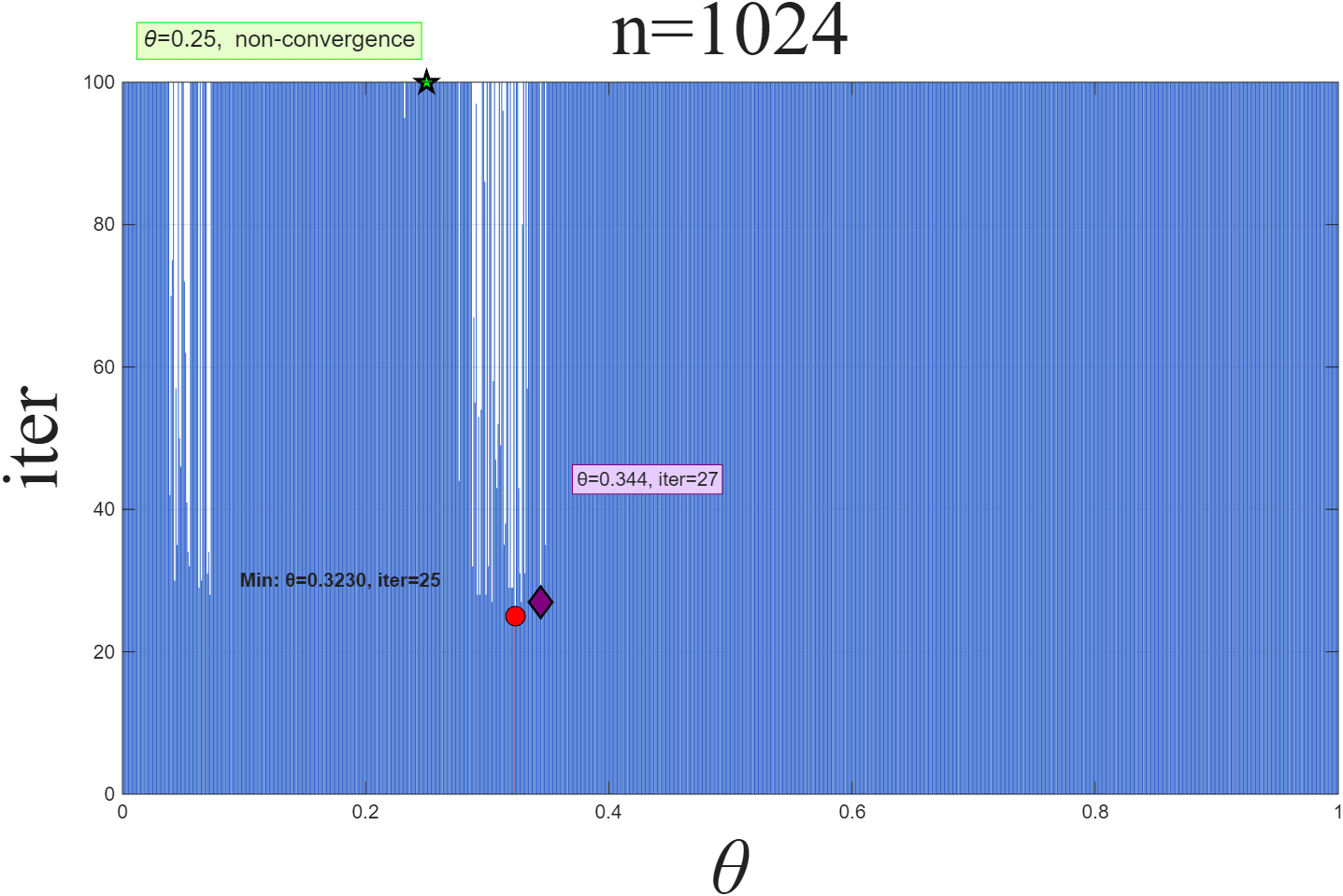}
	\end{minipage}
	\begin{minipage}{0.32\linewidth}
		\centering
		\includegraphics[width=1.1\linewidth]{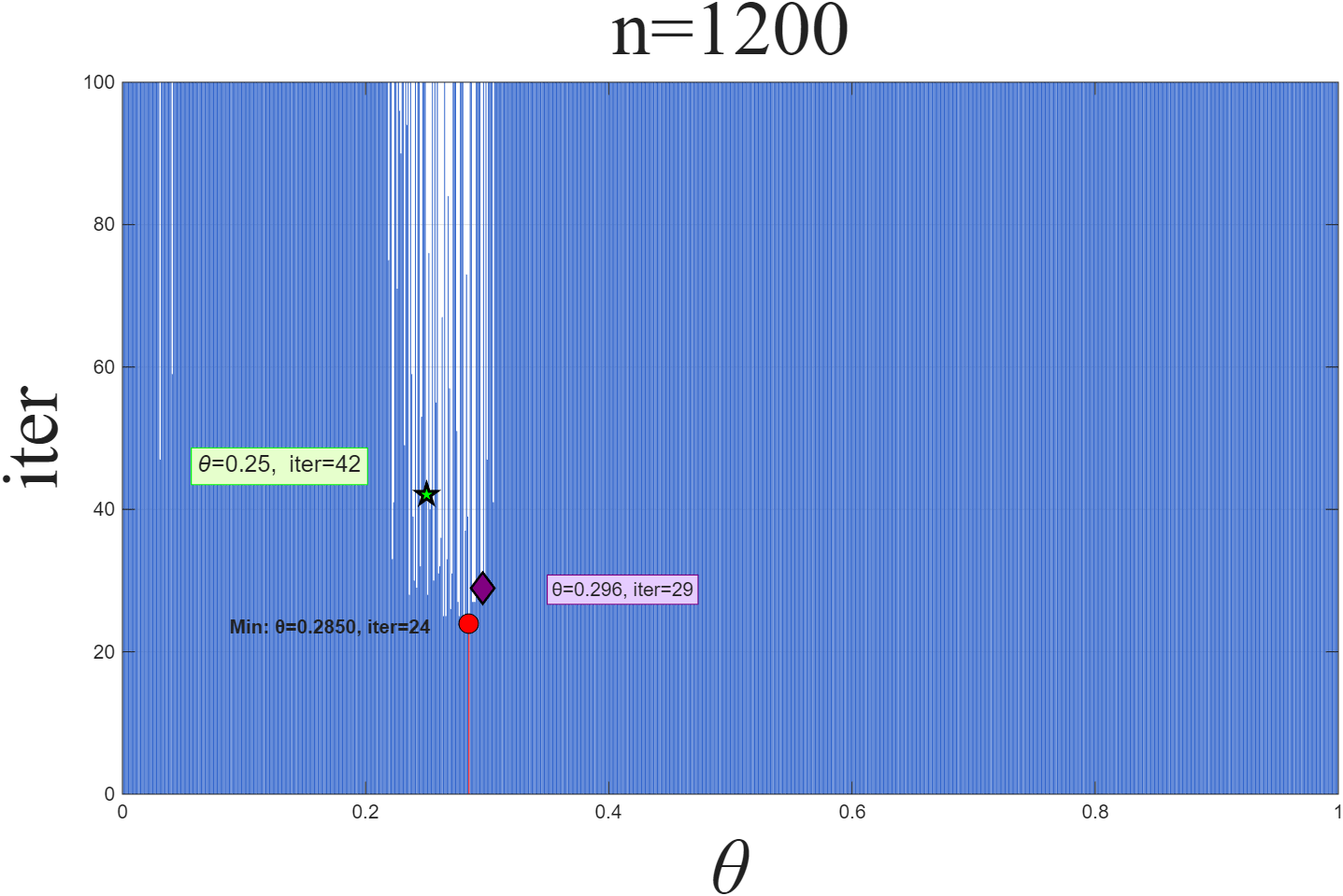}
	\end{minipage}
	\begin{minipage}{0.32\linewidth}
		\centering
		\includegraphics[width=1.1\linewidth]{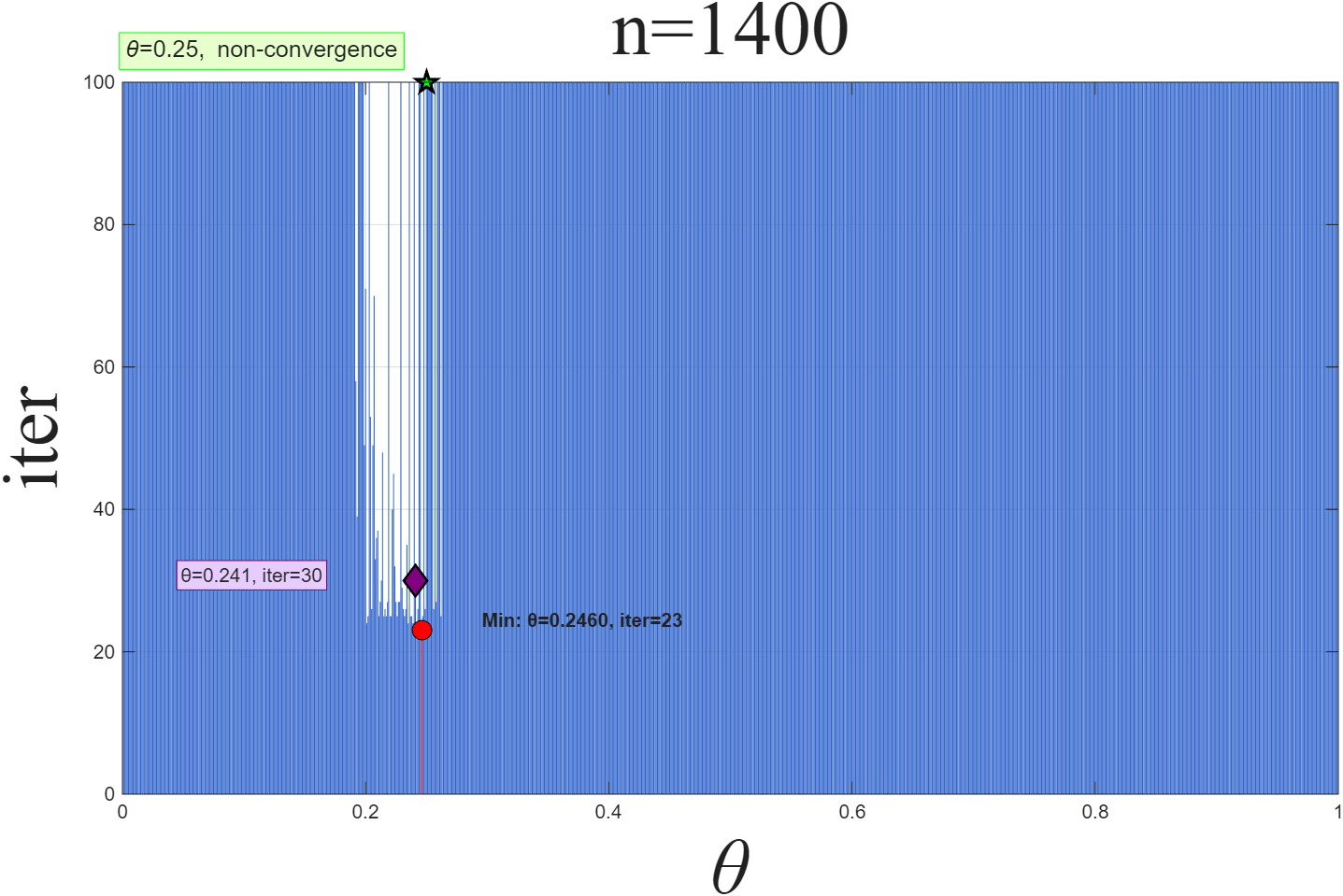}
	\end{minipage}
		\begin{minipage}{0.32\linewidth}
		\centering
		\includegraphics[width=1.1\linewidth]{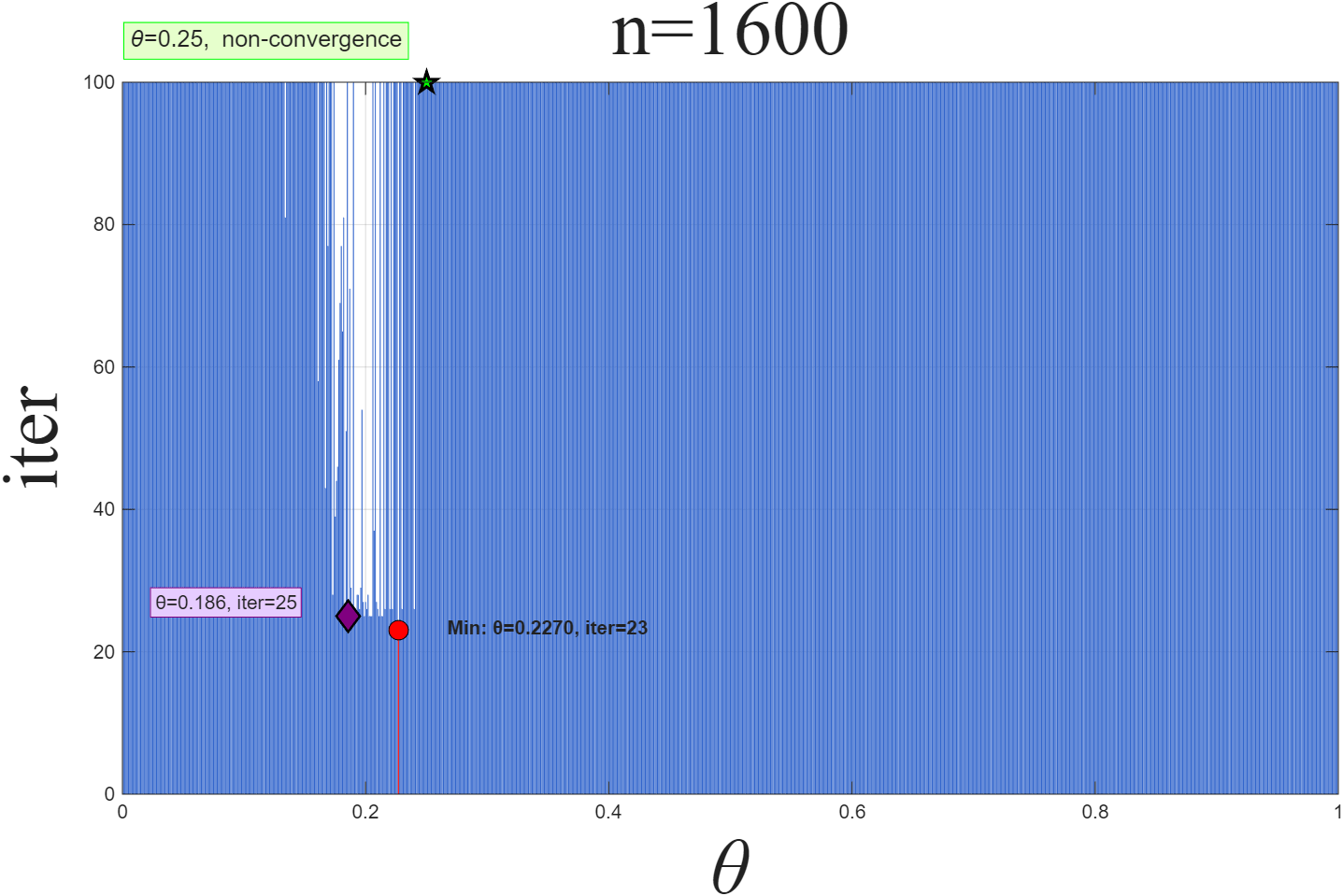}
	\end{minipage}
			\begin{minipage}{0.32\linewidth}
		\centering
		\includegraphics[width=1.1\linewidth]{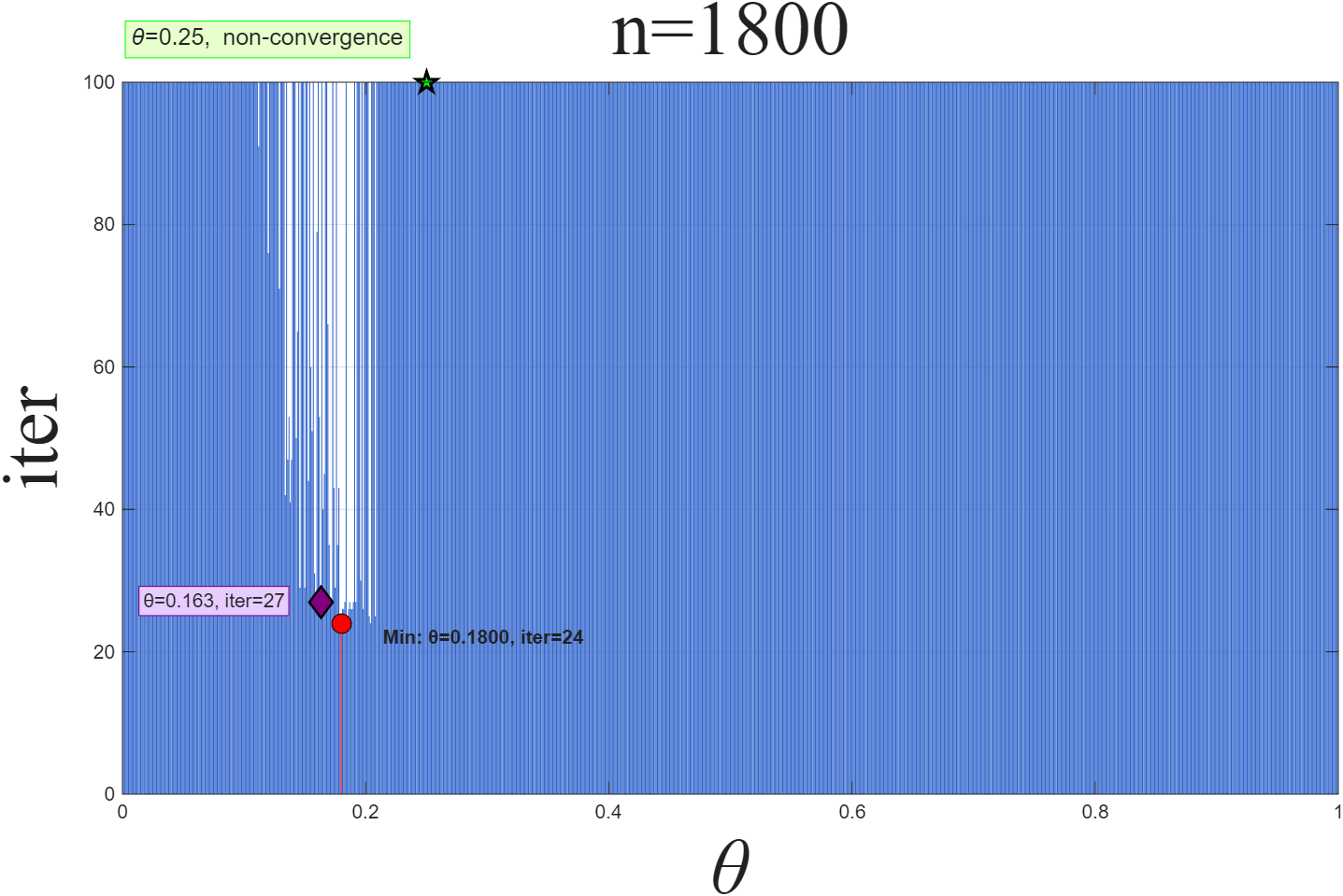}
	\end{minipage}
				\begin{minipage}{0.32\linewidth}
		\centering
		\includegraphics[width=1.1\linewidth]{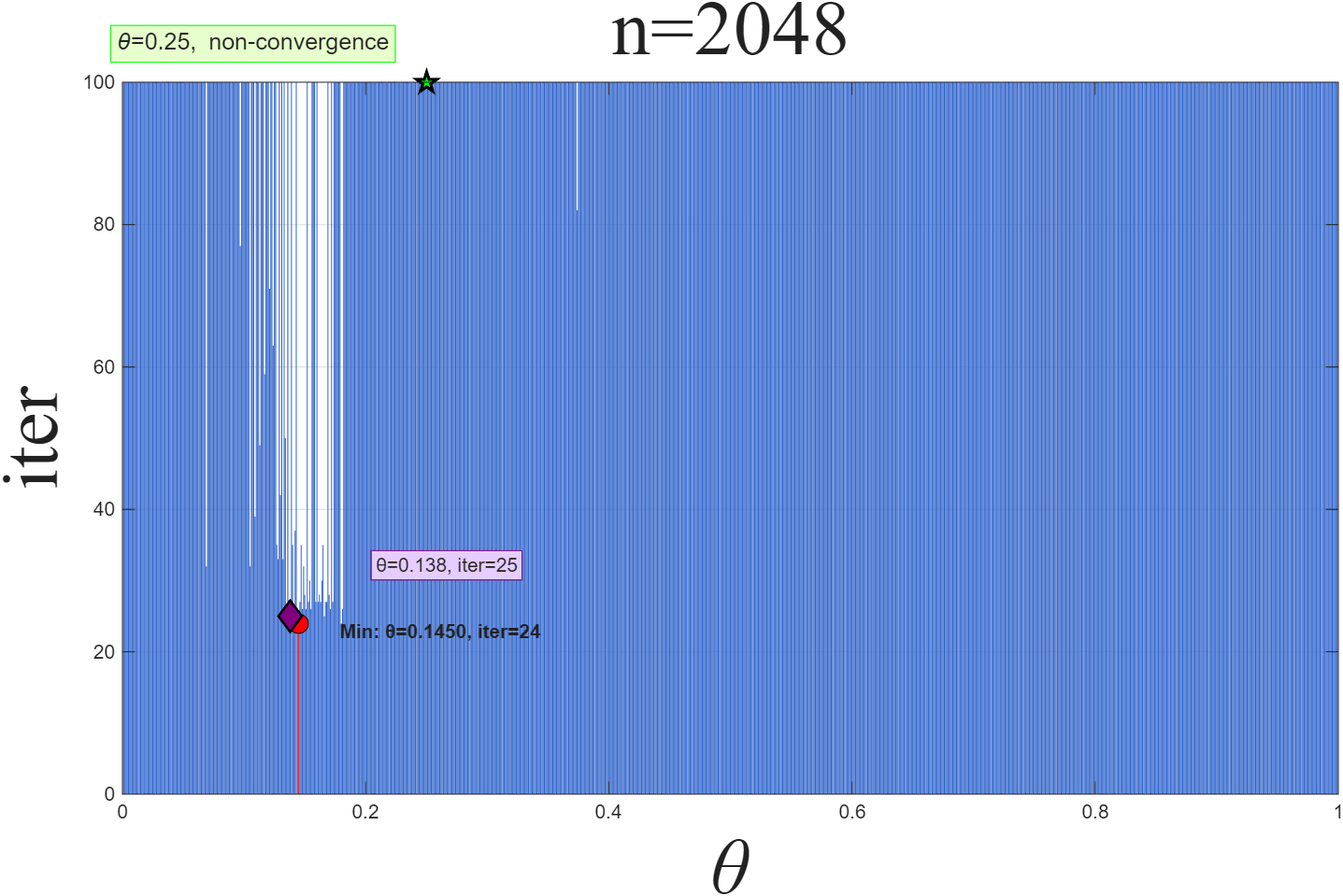}
	\end{minipage}
\end{figure}

Table \ref{hm_tab_pre} lists the result of the iteration counts and CPU time under three parameter selection schemes for the case of {  $n=1024,~ 1200,~ 1400,~ 1600,~1800,~2048$}, where optimal $\theta$ is the traversal-derived parameter. The parameters predicted by the GPR model exhibit high accuracy, differing from the optimal parameters by merely 0.021, while also closely matching the results obtained from the exhaustive search method.
Notice that AMG doesn't converge with default parameters in {  $n =1024,~ 1400,~1600,~1800,~2048$}, but it works well with predicted $\theta$ obtained by GPR. The number of iterations of AMG under the predicted parameter differs by only 2 from that under the optimal parameter and the CPU time increases by only 8.85$\%$. {  The traversal results under more diverse values of $\theta$ are detailed in Figure \ref{hm_travel_more_Helm}. As can be observed from the figure, for the Helmholtz equation, the parameter values of $\theta$ that lead to AMG convergence are concentrated within a very narrow region, which is accurately captured by GPR.}


\begin{table}[h]
	\centering
{ 
	\caption{ Comparison of results for AMG solving the Helmholtz equation with coefficient $2\pi$ }
	\label{hm_tab_pre}
	\scalebox{0.85}{ 
	\begin{tabular}{|c|c|c|c|c|c|c|c|c|c|}
		\hline
		$n$ & \makecell{predicted\\ $\theta$}  & iter &\makecell{CPU \\time} & \makecell{optimal\\ $\theta$} &iter & \makecell{CPU time} & \makecell{iter of \\default value} & \makecell{CPU time \\of default}   \\ \hline
			1024  &0.344&27   &{  4.5296} & 0.323 & 25 &{  4.3014} &\makecell{non-convergence} &\makecell{non-convergence} \\ \hline
		1200  &0.296&29   &{  6.5198} & 0.285 & 24 &{  6.4952} &42 &{  7.6307} \\ \hline
		1400  &0.241 &30 &{ 7.7401}  & 0.246 &23  &{ 7.7407} &\makecell{non-convergence} &\makecell{non-convergence} \\ \hline
		1600  &0.186   &25 &{  10.5838} & 0.2227 &23 &{  10.0305} & \makecell{non-convergence} &\makecell{non-convergence}  \\ \hline
	1800 &0.163   &27& {  14.1669} & 0.18 &24 &{  13.1320} & \makecell{non-convergence} &\makecell{non-convergence}  \\ \hline
		2048 &0.138  &25 &{  17.7594} & 0.145 &24 &{  17.2984} & \makecell{non-convergence} &\makecell{non-convergence} \\ \hline
	\end{tabular}
	}
	}
\end{table}

The Table \ref{Tab:hm:pre1} and \ref{Tab:hm:pre2} provide a comparative analysis of the predictive accuracy achieved by various combinations of kernel functions when applied to the Helmholtz equation with the coefficient {\( k = 2\pi \)}. 

{Table \ref{Tab:hm:pre1} presents the evaluation metrics of different kernel functions. As can be seen from the table, most dual-kernel combination functions exhibit higher $R^2$ and Corr values compared to the single Gaussian kernel function. It is particularly noteworthy that the Gaussian kernel achieves the lowest Corr value among all kernel functions. Specifically, the combination of Rational Quadratic and Laplacian kernels demonstrates optimal comprehensive performance: it achieves the lowest values in RMSE, MAE, BIC, and MdAPE, while simultaneously attaining the highest $R^2$ value. These results fully demonstrate that combined kernel functions can significantly enhance prediction performance compared to single kernel functions.}

Table \ref{Tab:hm:pre2} presents the PICP metrics for different kernel functions. The results show that most dual-kernel combined functions achieve higher PICP values than the single Gaussian kernel function. This further confirms that the kernel functions constructed through dual-kernel combinations can exhibit superior performance characteristics compared to single kernel function.

\begin{table}[h]
\centering
	\caption{ Evaluation of kernel function combinations for Helmholtz equation with coefficient $2\pi$}
	\label{Tab:hm:pre1}
	\scalebox{0.65}{
		\begin{tabular}
			{|c|c|c|c|c|c|c|c|c|}
			\hline \text { Kernel function } & \text { MSE } & \text { RMSE } & \text { MAE } & R$^{ 2 } $ & \text { BIC } & \text { Corr } & \text { MdAPE } & \text { LOO-SPE } \\ \hline
			\text { Gaussian+Laplacian } &5.9361 $\times 10^{-3}$ & 0.0770 & 0.0652 & 0.0461 & -150.6211 & 0.9271 & 0.3080 & 2.3634 $\times 10^{-3}$\\ \hline
			\text { Gaussian+Exponential } & 9.1618 $\times 10^{-3}$ & 0.0957 & 0.0802 & 0.0255 & -153.9052 & 0.7689 & 0.3347 & 2.8059 $\times 10^{-3}$ \\ \hline
			\text { Gaussian } & 3.9714 $\times 10^{-3} $ & 0.0630 & 0.0526 & 0.0224 & -251.7541 & 0.6933 & 0.2432 & 2.3984 $\times 10^{-4}$
			\\ \hline
			\text { Rational Quadratic+Laplacian } & 3.6738 $\times 10^{-3}$ & 0.0606 & 0.0442 & 0.0957 & -377.5221 & 0.7645 & 0.1425 & 2.2405 $\times 10^{-2}$ \\ \hline
			\text { Mat\"{e}rn+Laplacian } & 6.4141 $\times 10^{-3}$ & 0.0801 & 0.0673 & 0.5787 & -145.9441 & 0.9307 & 0.3139 & 2.4222 $\times 10^{-3} $
			\\ \hline
			\text { Rational Quadratic+Gaussian } & 4.6182 $\times 10^{-3} $& 0.0679 & 0.0612 & 0.0136 & -155.5238 & 0.9437 & 0.2517 & 2.2826 $\times 10^{-3}$ \\ \hline
			\text { Mat\"{e}rn+Gaussian+Laplacian } & 6.5222 $\times 10^{-3}$ & 0.0807 & 0.0679 & 0.0605 & -137.4683 & 0.9326 & 0.3152 & 2.4398 $\times 10^{-3} $\\ \hline
		\end{tabular}}
\end{table}

\begin{table}[h]
	\centering
	\caption{PICP evaluation for the Helmholtz equation with coefficient $2\pi$}
		\label{Tab:hm:pre2}
	\begin{tabular}{|c|c|}
		\hline
		Kernel function & PICP \\ \hline
		Gaussian+Laplacian & 85.7 \% \\ \hline
		Gaussian+Exponential & 85.7 \% \\ \hline
		Gaussian & 42.9 \% \\ \hline
		Rational Quadratic+Laplacian & 28.6\% \\ \hline
		Mat\"{e}rn+Laplacian & 71.4 \% \\ \hline
		Rational Quadratic+Gaussian & 85.7 \% \\ \hline
		Mat\"{e}rn+Gaussian+Laplacian & 85.7 \% \\ \hline
	\end{tabular}
\end{table}

\section{Summary}

In this paper, we have highlighted the effectiveness of kernel learning methods, particularly GPR, in optimizing parameters for AMG when solving various PDEs. The results demonstrated significant time savings and high reliability in parameter prediction, especially when using combined kernel functions. These findings provided valuable theoretical support for AMG optimization and practical applications.

In future, there are three promising directions worth exploring. First, mixed precision computing can be integrated with GPR-based AMG optimization. By leveraging low-precision computations for initial parameter exploration and high-precision adjustments for accuracy, computational efficiency can be further enhanced. Second, extending this approach to more complex PDEs, such as nonlinear equations and multiphysics problems, will require developing advanced kernel functions and multitask learning frameworks to handle their nonlinearity and multiscale nature. This extension could significantly improve the efficiency and accuracy of solving these challenging problems.
Third, we will further explore direct parameter prediction mechanisms based on PDE types, with the aim of fundamentally avoiding the need for repeated training procedures when encountering new PDEs. Besides, we will consider to employ Online GPR to address the challenges of few-shot learning and dynamic adaptation in large-scale matrix parameter prediction. This method utilizes an incremental learning mechanism to achieve dynamic optimization and real-time updating of model parameters, effectively overcoming the parameter adjustment latency inherent in traditional batch processing methods when handling streaming data.


In conclusion, the integration of kernel learning with AMG parameter optimization offers a promising path forward. However, a comprehensive theoretical explanation for the predictive efficacy of different kernels and the factors influencing prediction accuracy will be a key focus of future work. Furthermore, our subsequent research will focus on mixed precision  techniques on GPU architectures\cite{H-2022} ,   and their performance improvements for     AMG. Investigating complex PDEs with these advanced methods will likely yield more efficient and reliable numerical solutions for real-world scientific and engineering problems.

\section*{Acknowledgments}
{This work was supported in part by the High Performance Computing Platform of Xiangtan University.
The authors would like to thank the National Key R\&D Program of China (2023YFB3001604) for its support and assistance, and to thank the associate editor and the referees for their constructive comments and suggestions, which substantially improved the paper.}


\begin{thebibliography}{99}
	
	\bibitem{ref9} Antonietti P F, Caldana M, Dede L. Accelerating algebraic multigrid methods via artificial neural networks. Vietnam Journal of Mathematics, 51 (2023), pp. 1-36.
	
	\bibitem{ref15} Antonietti P F, Melas L. Algebraic multigrid schemes for high-order nodal discontinuous Galerkin methods. SIAM Journal on Scientific Computing, 42 (2020), pp. A1147-A1173.
	
	 \bibitem{CL2008} Chen L. iFEM: an innovative finite element methods package in MATLAB. Preprint, University of Maryland, 20, (2008).
	
	\bibitem{C-2006} Williams C K I,  Rasmussen C E. Gaussian processes for machine learning. Cambridge, MA: MIT press, (2006).
	
		
	\bibitem{F-Y-2022} Falgout R D, Yang U M. hypre: A library of high performance preconditioners. International Conference on computational science. Berlin, Heidelberg: Springer Berlin Heidelberg, (2002),pp. 632-641.
	
	\bibitem{ref25} Gronau Q F, Wagenmakers E J. Limitations of Bayesian leave-one-out cross-validation for model selection. Computational Brain \& Behavior, 2 (2019), pp. 1-11.
	
	\bibitem{ref12} Gottschalk H, Kahl K. Coarsening in algebraic multigrid using Gaussian processes. Electronic Transactions on Numerical Analysis, 54 (2021), pp. 514-533.
	
	\bibitem{ref10} Haghi P, Geng T, Guo A, Wang T, Herbordt M. FP-AMG: FPGA-based acceleration framework for algebraic multigrid solvers. In: 2020 IEEE 28th Annual International Symposium on Field-Programmable Custom Computing Machines (FCCM), (2020), pp. 148-156.
	
	\bibitem{ref24} Hodson T O. Root mean square error (RMSE) or mean absolute error (MAE): When to use them or not. Geoscientific Model Development Discussions, (2022), pp. 1-10.
	
	\bibitem{ref16} Jiang K, Su X, Zhang J. A general alternating-direction implicit framework with Gaussian process regression parameter prediction for large sparse linear systems. SIAM Journal on Scientific Computing, 44 (2022), pp. A1960-A1988.
	
	\bibitem{ref17} Jiang K, Zhang J, Zhou Q. Multitask kernel-learning parameter prediction method for solving time-dependent linear systems. CSIAM Transactions on Applied Mathematics, 4 (2023), pp. 672-695.
	
	\bibitem{ref3} Kickinger F. Algebraic multi-grid for discrete elliptic second-order problems. In: Multigrid Methods V: Proceedings of the Fifth European Multigrid Conference held in Stuttgart, (1998), pp. 157-172.
	
	\bibitem{ref1} Ruge J W, St\"{u}ben K. Algebraic multigrid. In: Multigrid Methods, (1987), pp. 73-130.
	
	\bibitem{ref4} Richter C, Schps S, Clemens M. GPU acceleration of algebraic multigrid preconditioners for discrete elliptic field problems. IEEE Transactions on Magnetics, 50 (2014), pp. 461-464.
	

	\bibitem{ref8} Luz I, Galun M, Maron H, Basri R, Yavneh I. Learning algebraic multigrid using graph neural networks. In: International Conference on Machine Learning, (2020), pp. 6489-6499.
	
	\bibitem{ref5} Brezina M, Falgout R, MacLachlan S, Mccormick T, Mccormick S, Ruge J. Adaptive algebraic multigrid. SIAM Journal on Scientific Computing, 27 (2006), pp. 1261-1286.
	
	\bibitem{ref26} Neath A A, Cavanaugh J E. The Bayesian information criterion: background, derivation, and applications. Wiley Interdisciplinary Reviews: Computational Statistics, 4 (2012), pp. 199-203.
	
	\bibitem{ref2} Bastian P, Blatt M, Scheichl R. Algebraic multigrid for discontinuous Galerkin discretizations of heterogeneous elliptic problems. Numerical Linear Algebra with Applications, 19 (2012), pp. 367-388.
	
	\bibitem{ref23} Petit S J, Bect J, Feliot P, Vazquez E. Parameter selection in Gaussian process interpolation: an empirical study of selection criteria. SIAM/ASA Journal on Uncertainty Quantification, 11 (2023), pp. 1308-1328.
	
	\bibitem{ref27} Rad K R, Maleki A. A scalable estimate of the out-of-sample prediction error via approximate leave-one-out cross-validation. Journal of the Royal Statistical Society Series B: Statistical Methodology, 82 (2020), pp. 965-996.

	\bibitem{ref19} Rahmati O, Choubin B, Fathabadi A, Coulon F, Soltani E, Shahabi H, Mollaefar E, Tiefenbacher J, Cipullo S, Ahmad BB.  Predicting uncertainty of machine learning models for modelling nitrate pollution of groundwater using quantile regression and UNEEC methods. Science of the Total Environment, 688 (2019), pp. 855-866.
	
	\bibitem{ref21} Ranjan R, Huang B, Fatehi A. Robust Gaussian process modeling using EM algorithm. Journal of Process Control, 42 (2016), pp. 125-136.
	
	\bibitem{ref14} Sbai M A, Larabi A. On solving groundwater flow and transport models with algebraic multigrid preconditioning. Groundwater, 59 (2021), pp. 100-108.
	
	\bibitem{ref30} Schulz E, Speekenbrink M, Krause A. A tutorial on Gaussian process regression: Modelling, exploring, and exploiting functions. Journal of Mathematical Psychology, 85 (2018), pp. 1-16.
	
	\bibitem{ref20} Shrestha D L, Solomatine D P. Machine learning approaches for estimation of prediction interval for the model output. Neural Networks, 19 (2006), pp. 225-235.
	
	\bibitem{ref7} Taghibakhshi A, MacLachlan S, Olson L, West M. Optimization-based algebraic multigrid coarsening using reinforcement learning. Advances in Neural Information Processing Systems, 34 (2021), pp. 12129-12140.
	
	\bibitem{ref28} Vehtari A, Gelman A, Gabry J. Practical Bayesian model evaluation using leave-one-out cross-validation and WAIC. Statistics and Computing, 27 (2017), pp. 1413-1432.
	
	\bibitem{ref29} Watanabe S. A widely applicable Bayesian information criterion. The Journal of Machine Learning Research, 14 (2013), pp. 867-897.
	
	\bibitem{ref11} Wiesner T A, Mayr M, Popp A, Gee M W, Wall W A. Algebraic multigrid methods for saddle point systems arising from mortar contact formulations. International Journal for Numerical Methods in Engineering, 122 (2021), pp. 3749-3779.
	
	\bibitem{ref18} Wan X, Li X, Wang X, Yi X, Zhao Y, He X, Wu R, Huang M.  Water quality prediction model using Gaussian process regression based on deep learning for carbon neutrality in papermaking wastewater treatment system. Environmental Research, 211 (2022), pp. 112942.
	
	\bibitem{ref6} Notay Y. An aggregation-based algebraic multigrid method. Electronic Transactions on Numerical Analysis, 37 (2010), pp. 123-146.
	
	\bibitem{ref13} Ye S, Xu X, An H, Yang X. A supplementary strategy for coarsening in algebraic multigrid. Applied Mathematics and Computation, 394 (2021), pp. 125795.
	
	\bibitem{ref22} Yadav A, Bareth R, Kochar M, Pazoki M, Sehiemy R A E. Gaussian process regression-based load forecasting model. IET Generation, Transmission \& Distribution, 18 (2024), pp. 899-910.
	
	\bibitem{Z-X-2024}  Zou H, Xu X, Zhang C S, Mo Z Y.
	AutoAMG ($\theta $): An Auto-tuned AMG Method Based on Deep Learning for Strong Threshold. Communications in Computational Physics, 36 (2024), pp. 200-220. 
	\bibitem{D-2008}{     de Sterck H, Falgout R D, Nolting J W, Yang, U M. Distance‐two interpolation for parallel algebraic multigrid. Numerical Linear Algebra with Applications, 15 (2008), pp. 115-139.          }






\bibitem{H-2022} { 
Higham, N J,  Mary, T. Mixed precision algorithms in numerical linear algebra. Acta Numerica, 31 (2022), pp. 347-414.}
	
	
\end{thebibliography}
\end{document}